\documentclass[11pt,leqno]{article} 
\begin{document}
%
\baselineskip=1.1\baselineskip 
%
\newcommand{\cC}{{\cal C}}
\newcommand{\cB}{{\cal B}}
\newcommand{\cE}{{\cal E}}
\newcommand{\cH}{{\cal H}}
\newcommand{\cK}{{\cal K}}
\newcommand{\cL}{{\cal L}}
\newcommand{\cI}{{\cal I}}
\newcommand{\cM}{{\cal M}}
\newcommand{\cN}{{\cal N}}
\newcommand{\cF}{{\cal F}}
\newcommand{\cG}{{\cal G}}
\newcommand{\cP}{{\cal P}}
\newcommand{\cR}{{\cal R}}
\newcommand{\cS}{{\cal S}}
\newcommand{\cT}{{\cal T}}
\newcommand{\cV}{{\cal V}}
\newcommand{\TT}{{\rm Tr}}
\newcommand{\ds}{\displaystyle}
\newcommand{\sss}{\scriptstyle}
\newcommand{\rr}{I\!\!R}
\newcommand{\D}{I\!\!D}
\newcommand{\dd}{D^{\mu}}
\newcommand{\nn}{I\!\!N}
\newcommand{\cc}{I\!\!C}
\newcommand{\LL}{I\!\!L}
\newcommand{\un}{\underline}
\newcommand{\hs}{\hspace}
\newcommand{\vs}{\vspace}
\newcommand{\no}{\noindent}
\newcommand{\be}{\begin{equation}}
\newcommand{\ee}{\end{equation}}
\newcommand{\ba}{\begin{array}}
\newcommand{\ena}{\end{array}}
\newcommand{\diag}{\mbox{\rm diag}}
\newcommand{\1}{\mbox{\bf 1}}
\newcommand{\bdes}{\begin{description}}
\newcommand{\edes}{\end{description}}
\newcommand{\bla}{\begin{lemma}}
\newcommand{\ela}{\end{lemma}}
\newcommand{\bcy}{\begin{corollary}}
\newcommand{\ecy}{\end{corollary}}
\newtheorem{thm}{Theorem}[section]   
\newtheorem{re}[thm]{Remark}
\newtheorem{co}[thm]{Corollary}
\newtheorem{pr}[thm]{Proposition}
\newtheorem{de}[thm]{Definition}
\newtheorem{lm}[thm]{Lemma}
\newtheorem{example}[thm]{Example}
\renewcommand{\theequation}{\thesection.\arabic{equation}}
\def\squarebox#1{\hbox to #1{\hfill\vbox to #1{\vfill}}}
\newcommand{\qed}{\hspace*{\fill}
\vbox{\hrule\hbox{\vrule\squarebox{.45em}\vrule}\hrule}\smallskip}
\def \ci {\mathop{\hbox {\vrule height .4pt depth 0pt width 0.5cm
\hskip -0.3cm \vrule height 10pt
depth 0pt \hskip 0.1cm \vrule height 10pt depth 0pt \hskip 0.2cm}}}
\def\E{\mbox{\rm E}}
\def\supp{\mbox{\rm supp}}
%
\newdimen\fletxauu
\newdimen\fletxadosu
\newdimen\fletxaudos
\newdimen\fletxadosdos
\newdimen\difprimera
\newdimen\difsegona
\newdimen\diftercera
\newbox\Zero
\newbox\Una
\newbox\Dues
\newbox\Tres
\newbox\Quatre
\newbox\Cinc
\newbox\Sis
\newbox\Set
\newbox\Vuit
\newbox\Nou
\def\mapstofill{$\mathsurround=0pt\mapstochar\mathrel{\mkern-4mu}
                \mathord- \mkern-6mu
                \cleaders\hbox{$\mkern-2mu\mathord-\mkern-2mu$}
                \hfill\mkern-6mu\mathord\rightarrow$}
\def\composiciodeaplicacions: nom: #1 conjunt1: #2 aplicacio1: #3
     conjunt2: #4 aplicacio2: #5 conjunt3: #6 element1: #7
     element2: #8 element3: #9 fi{
    \fletxauu=1.5cm
    \fletxadosu=1.5cm
    \fletxaudos=1.5cm
    \fletxadosdos=1.5cm
    \setbox\Una=\hbox{#1}
    \setbox\Dues=\hbox{#2}
    \setbox\Tres=\hbox{#3}
    \setbox\Quatre=\hbox{#4}
    \setbox\Cinc=\hbox{#5}
    \setbox\Sis=\hbox{#6}
    \setbox\Set=\hbox{#7}
    \setbox\Vuit=\hbox{#8}
    \setbox\Nou=\hbox{#9}
    \difprimera=\wd\Tres
    \advance\difprimera by -0.4cm
    \ifdim\difprimera>0pt\advance\fletxauu by \difprimera
                         \advance\fletxadosu by \difprimera
                         \fi
    \difprimera=\wd\Dues
    \advance\difprimera by -\wd\Set
    \divide\difprimera by 2
    \ifdim\difprimera>0pt\advance\fletxadosu by \difprimera
                  \else\advance\fletxauu by -\difprimera
          \fi
    \difsegona=\wd\Cinc
    \advance\difsegona by -0.4cm
    \ifdim\difsegona>0pt\advance\fletxaudos by \difsegona
                        \advance\fletxadosdos by \difsegona
                        \fi
    \difsegona=\wd\Quatre
    \advance\difsegona by -\wd\Vuit
    \divide\difsegona by 2
    \ifdim\difsegona>0pt\advance\fletxadosu by \difsegona
                      \advance\fletxadosdos by \difsegona
                 \else\advance\fletxauu by -\difsegona
                      \advance\fletxaudos by -\difsegona
          \fi
    \diftercera=\wd\Sis
    \advance\diftercera by -\wd\Nou
    \divide\diftercera by 2
    \ifdim\diftercera>0pt\advance\fletxadosdos by \diftercera
                  \else\advance\fletxaudos by -\diftercera
          \fi
   %
$$
    \begin{array}{l}
    \copy\Dues\vbox{\hsize=\fletxauu
                     \offinterlineskip
                     \hbox to \fletxauu{\hfil\copy\Tres\hfil}
                     \kern0pt
                     \hbox to \fletxauu{\kern1.5pt
                                        \rightarrowfill
                                        \kern1.5pt}}
      \copy\Quatre\vbox{\hsize=\fletxaudos
                      \offinterlineskip
                      \hbox to \fletxaudos{\hfil\copy\Cinc\hfil}
                      \kern0pt
                      \hbox to \fletxaudos{\kern1.5pt
                                           \rightarrowfill
                                           \kern1.5pt}}
      \copy\Sis\\ [2mm]
    \hskip\difprimera
      \copy\Set\hbox to \fletxadosu{\kern1.5pt
                                    \mapstofill
                                    \kern1.5pt}
      \copy\Vuit\hbox to \fletxadosdos{\kern1.5pt
                                      \mapstofill
                                      \kern1.5pt}
      \copy\Nou
    \end{array}
$$
    }
\def\aplicacio: nom: #1 conjunt1: #2 aplicacio: #3
     conjunt2: #4 element1: #5
     element2: #6 fi{
    \fletxauu=1.5cm
    \fletxadosu=1.5cm
    \setbox\Una=\hbox{#1}
    \setbox\Dues=\hbox{#2}
    \setbox\Tres=\hbox{#3}
    \setbox\Quatre=\hbox{#4}
    \setbox\Cinc=\hbox{#5}
    \setbox\Sis=\hbox{#6}
    \difprimera=\wd\Tres
    \advance\difprimera by-0.4cm
    \ifdim\difprimera>0pt\advance\fletxauu by \difprimera
                         \advance\fletxadosu by \difprimera
                         \fi
    \difprimera=\wd\Dues
    \advance\difprimera by -\wd\Cinc
    \divide\difprimera by 2
    \ifdim\difprimera>0pt\advance\fletxadosu by \difprimera
                    \else\advance\fletxauu by -\difprimera
          \fi
    \difsegona=\wd\Quatre
    \advance\difsegona by -\wd\Sis
    \divide\difsegona by 2
    \ifdim\difsegona>0pt\advance\fletxadosu by \difsegona
                   \else\advance\fletxauu by -\difsegona
          \fi
$$
    \begin{array}{l}
    \ifdim\wd\Una>0pt\copy\Una \colon \,\,
              \else\fi
     \copy\Dues \vbox{\hsize=\fletxauu
                             \offinterlineskip
                             \hbox to \fletxauu{\hfil \copy\Tres \hfil}
                             \kern0pt
                             \hbox to \fletxauu{\kern1.5pt
                                                \rightarrowfill
                                                \kern1.5pt}}
                \copy\Quatre \\ [2mm]
                \hskip\difprimera
                \copy\Cinc \hbox to \fletxadosu{\kern1.5pt
                                            \mapstofill
                                            \kern1.5pt}
                \copy\Sis
    \end{array}
$$
    }
\title{Linear stochastic differential equations
\\
with functional boundary conditions\vs{1cm}}
\author{
Aureli Alabert
\thanks{Supported by grants SGR99-87 of CIRIT and BFM2000-0009 of DGESIC}
\\
Departament de Matem\`atiques \\
Universitat Aut\`onoma de Barcelona \\
08193 Bellaterra, Catalonia \\
e-mail: alabert@mat.uab.es
\and
Marco Ferrante
\thanks{Supported by grant COFIN9901244421 of MURST}
\\
Dip. di Matematica Pura ed Appl. \\
Universit\`a degli Studi di Padova \\
via Belzoni 7, 
35131 Padova, Italy \\
e-mail: ferrante@math.unipd.it}
\maketitle
\begin{abstract}
\no
We consider linear $n$-th order
stochastic differential equations on $ [0,1] $,
with linear boundary conditions supported by a finite subset of $[0,1]$.
We study some features of the 
solution to these problems, 
and especially its 
conditional independence properties of Markovian type.
\end{abstract}

\vspace{3truecm}
  {\bf AMS Classification:} 60H10, 60J25

\thispagestyle{empty}
\vfill\eject
\null
\thispagestyle{empty}
\vfill\eject
\setcounter{page}{1}
\section{Introduction}
  It is well known that, under suitable Lipschitz and growth 
  conditions on the coefficients, 
  a classical It\^{o} stochastic differential equation
\begin{equation}\label{a-itoeq}
  X(t)=\xi+\int_0^t b(s,X(s))\,ds+\int_0^t \sigma(s,X(s))\,dW(s)
  \ ,
\end{equation}
  where $W$ is a Wiener process and $\xi$ is a ${\cal F}_0$-measurable
  random variable for a given non-anticipating filtration
  $\{{\cal F}_t,\ t\ge 0\}$ of $W$, 
  has a unique strong solution
  which is a Markov process.
\par
  If $\xi$ is not ${\cal F}_0$-measurable or the coefficients $b$, $\sigma$
  are random and non-adapted, then any
  reasonable interpretation of $X$ in (\ref{a-itoeq}) will not be an 
  ${\cal F}_t$-adapted process and, unless $\sigma$ is a constant, 
  we need to use some anticipating
  stochastic integral to give a sense to the equation.
  In these cases, the solution is not a Markov process in
  general.
\par
  Still another setting that leads to anticipation is the case
  of boundary conditions. That means, the first variable of the
  solution process is no longer a datum of the problem, 
  time runs in a bounded
  interval, say from 0 to 1, and we impose a relation
  $h(X(0),X(1))=0$ between the first and the last variables of
  the solution. In this situation, the fact that the solution will not
  be Markovian is quite intuitive, since the strong relationship
  between $X(0)$ and $X(1)$ will prevent the independence
  of $X(0)$ and $X(1)$ from holding, even when conditioning to $X(a)$,
  $a\in\mathopen]0,1\mathclose[$, except maybe in some very particular cases.
\par
  On the other hand, it may also seem
  intuitive that
  the following weaker conditional independence property can hold true:
  For any $0\le a<b\le 1$, the $\sigma$-fields 
  $\sigma\{X(t),\ t\in[a,b]\}$
  and 
  $\sigma\{X(t),\ t\in\mathopen]a,b\mathclose[^c\}$
  are conditionally independent given $\sigma\{X(a),X(b)\}$. 
  We will denote it
by 
\begin{equation}\label{a-intuitci}
  \sigma\{X(t),\ t\in[a,b]\}
  \ci_{\sigma\{X(a),X(b)\}}
  \sigma\{X(t),\ t\in\mathopen]a,b\mathclose[^c\}
  \ .
\end{equation}
  Now $X(0)$ and $X(1)$ are on the same side in relation
  (\ref{a-intuitci}),
  so that the boundary condition does not seem to cause
  the problem seen above.
  But the following example shows that this is wrong:
\begin{example}{\label{a-NP1}}
  {\rm
  Consider the problem
\[
\left\{
\begin{array}{l}
{\ds
\dot X(t)=f(X(t))\,dt+\dot W(t)
\ ,\quad
t \in [0,1]
}
\\ [2mm]
 h(X(0),X(1))=0
 \ ,
\end{array}
\right.
\]
  where the noise appears additively, and assume
  that a unique solution exists and that the boundary
  condition given by $h$ does not
  reduce to an initial or final condition. Then, relation
  (\ref{a-intuitci}) holds if and only if $f(x)=\alpha x+\beta$, for some
  constants $\alpha$ and $\beta$. This was proved in Nualart and
  Pardoux \cite{np1}.
  \qed
  }
\end{example}
  The processes satisfying (\ref{a-intuitci}) were 
  called {\em reciprocal processes\/} 
  by S. Bernstein \cite{ber}. 
  The concept arose directly from E. Schr\"odinger ideas on the
  formulation of quantum mechanics.
  More recent research on such 
  processes has been carried out by B. Jamison \cite{jam},
  A. Krener \cite{kre}, R. Frezza, A. Krener and B. Levy
  \cite{FKL}, M. Thieullen \cite{thi} and J.C. Zambrini \cite{zam}.
\par
  Other names can be found in the literature to refer to the
  same concept. A reciprocal process is a {\sl one-parameter Markov
  field\/} in Paul L\'evy's terminology, and is also called
  a {\sl quasi-Markov process}, a {\sl local Markov process\/} and a
  {\sl Bernstein process}.
  We shall simply call them Markov fields (see Definition \ref{3a}).
\begin{example}\label{a-NP2}
  {\rm
  Consider now the problem
\[
\left\{
\begin{array}{l}
{\ds
\ddot X(t)+f(X(t),\dot X(t))=\dot W(t)
\ ,\quad
t \in [0,1]
}
 \\ [2mm]
 X(0)=c_1\ ,\ X(1)=c_2
 \ .
\end{array}
\right.
\]
  This is a second order stochastic differential equation, and
  it is natural to ask for conditional independence
  properties of the 2-dimensional process $Y(t)=(\dot X(t),X(t))$,
  since $X(t)$ has $C^1$ paths, and therefore it is meaningless to
  look for this kind of properties for $X(t)$ itself.
\par
  Nualart and Pardoux \cite{np2} proved that if $Y(t)$ is a Markov field,
  then, as in Example \ref{a-NP1}, $f$ must be an affine
  function. Moreover, if $f$ is affine, then $Y$ is not only
  a Markov field, but a Markov process.
  \qed
  }
\end{example}
  Let us look at this example more closely: Note that a
  boundary condition for a second order equation has the general form
$$
  h(Y(0),Y(1))=(0,0)\ ,
$$
  where $Y(0)=(\dot X(0),X(0))$ and $Y(1)=(\dot X(1),X(1))$. 
  However, in Example \ref{a-NP2} the two scalar
  conditions do not mix values at 0 and values at 1 of $Y$.
  The same happens, for instance, with the Neumann-type conditions
  $\dot X(0)=c_1$, $\dot X(1)=c_2$, and the result is the same
  ($Y$ Markov field $\Rightarrow$ $f$ affine $\Rightarrow$
   $Y$ Markov process).
\par
  From these examples and other equations of first and second order
  that have been studied so far
  (see e.g. \cite{op1}, \cite{afn}, \cite{an1}, \cite{am1}),
  we learn that
\par
\begin{enumerate}
\itemsep=-3pt           
\item The Markovian properties can be expected only in
  ``linear'' cases.
\item
  The specific Markovian property depends on the actual 
  form of the boundary condition.
\end{enumerate}
\par
  It should also be noted that the requirement of linearity on the drift
  coefficient $f$ is related to the fact that the noise appears
  additively. Should not this be the case, the Markovian 
  property would occur under a different condition
  which relates the drift and the diffusion coefficients
  (see \cite{afn} and \cite{am1}).
\par
\bigskip
  In the present paper we will consider linear stochastic differential
  equations of arbitrary order with additive white noise. Our boundary
  conditions will not be restricted to involve the solution process
  at the endpoints of the time interval, but we will allow them to involve
  the values at finitely many points inside the interval.
  They are usually called {\sl functional\/} or {\sl lateral\/} 
  boundary conditions. 
  Our main goal is to seek which kind of conditional independence
  properties can be established for the solution. A preliminary work in this
  direction was published in Alabert and Ferrante \cite{af1}. 
  Here we considerably refine and extend
  the results therein. 
  This type of equations was already
  considered by Russek \cite{rus}, who proved that the solutions are
  Markov processes if and only if the lateral conditions fix to a constant
  the variables $X(t)$, for all points $t$ in the support of the conditions.
  His techniques, based in the notion of reproducing kernel space,
  are different from ours.
\par
  Our main result (Theorem \ref{mainth}) can be stated in the following
  way: Fix two points $0\le a< b\le 1$ and set 
  $Y(t)=(D^{n-1}X(t),\dots,DX(t), X(t))$, where $n$ is the order of the equation, 
  $X(t)$ is its solution process, and $D$ is the time derivative;
  the process $\{Y(t),\ t\in [0,1]\}$ 
  satisfies the relation (\ref{a-intuitci})
  if and only if there are no lateral
  conditions involving points inside and outside the interval $[a,b]$.
  We also state a conditional
  independence property for the case when there are conditions that do involve
  points inside and outside $[a,b]$ (Theorem \ref{main3}). Finally we obtain 
  a result from which Russek's theorem can be trivially recovered
  (Theorem \ref{main2}). 
  The paper is organised as follows: 
\par
  In Section \ref{lsdewfbc} we precise 
  the statement of the problem and develop some notation and properties
  that will be needed later. 
\par
  Section 3 contains the main probabilistic tools: Lemma \ref{2g} and Proposition
  \ref{gogknud}. The first is 
  a characterisation of the conditional independence of two random vectors
  given a function of them. It is the most important ingredient in 
  the proof of Theorem \ref{mainth}, but it cannot be applied for certain 
  singular values of $a$ and $b$. 
  For these values, we employ an approximation by the solution of  
  perturbed equations. The approximation argument involves the convergence in 
  $L^2$ of 
  a sequence of conditional expectations with varying conditioning $\sigma$-fields.
  Proposition \ref{gogknud}
  gives a sufficient condition for this convergence in a general setting.
\par
  In Section 4 we establish the main results. The proofs of Theorems
  \ref{main3} and \ref{main2} will be only sketched, since the procedure is
  similar to that of Theorem \ref{mainth}, with slight modifications.

\setcounter{equation}{0}
\section{Linear SDE with functional boundary conditions}
\label{lsdewfbc}
The present section will be devoted to the statement of the problem,
the definition of a solution, 
and to absolute continuity and approximation results
for the solution of an
$ n $--th order linear stochastic differential equation
with linear functional boundary conditions.
\subsection{Statement of the problem and definition of a solution}
\par
Consider the differential operator
\[
L := D^n + a_{n-1} D^{n-1} + \cdots
+ a_{1} D + a_{0}
\ ,\quad
D :=  \frac{d}{dt}
\ ,
\]
where $ a_i $ are continuous functions on $[0,1]$.
Let $ \{ W(t),\ t \in [0,1] \}$ be
a standard Wiener process. We assume that $W$ is the coordinate
process in the classical Wiener space
$(\Omega,{\cal F},P)$, that means, 
$\Omega=C_0([0,1];\rr)$ is the space of continuous functions
on $[0,1]$ vanishing at zero, ${\cal F}$ its Borel $\sigma$-field, 
and $P$ the Wiener measure.
We shall deal with
the SDE
\be
\label{1}
L [X] = \dot{W}
\ee
on $[0,1]$,
together with the additional conditions
\be
\label{2}
\sum_{j=1}^{m} \alpha_{ij} X(t_j) = c_i
\ ,\quad
 1 \leq i \leq n
\ ,
\ee
where $ m \geq n $, $ 0 \leq t_1 < \ldots < t_m \leq 1 $ are some
given points in $[0,1]$,
and $ \alpha_{ij} $, $ c_i $ are real numbers.
The matrix of coefficients $(\alpha_{ij})$ is assumed to have
full rank.

As in the case of ordinary differential equations,
(\ref{1})-(\ref{2}) can be regarded as a first order system
\be
\label{4}
D Y(t) + A(t) Y(t) =
\dot B (t)
\ ,\quad
t \in [0,1]
\ ,
\ee
with constraints
\be
\label{a-fbcY}
\sum_{j=1}^{m} \alpha_{ij} Y_n(t_j) = c_i
\ ,\quad
 1 \leq i \leq n
\ ,
\ee
where
$ Y(t) = \big(Y_1(t), \ldots, Y_n(t) \big) $,
$ Y_i(t) = D^{n-i} X(t) $ for $ 1 \leq i \leq n $,
$B(t)=(W(t),0,\dots,0)$, and
\be
\label{5}
A(t) =
\left[
\begin{array}{cccccc}
a_{n-1} (t) &  a_{n-2} (t) & \cdots  & a_{1} (t) & a_{0} (t)
\\
-1          &  0           & \cdots  & 0         & 0
\\
0           &  -1          & \cdots  & 0         & 0
\\
\vdots      & \vdots       & \ddots  & \vdots    & \vdots
\\
0           & 0            & \cdots  & -1        & 0
\end{array}
\right]
\ .
\ee
The lateral condition
(\ref{a-fbcY}) is a special case of the general
linear condition
\be
\label{3}
\Lambda [ Y ] = c
\ ,
\ee
for 
$ \Lambda$ in the set ${\cal L}\left(C \left([0,1];\rr^n \right); \rr^n \right) $
of linear continuous $\rr^n$-valued functionals
on $ C \left( [0,1];\rr^n \right) $, and
$ c \in \rr^n $.
By the Riesz representation
theorem, (\ref{3}) can be written as    
\be
\label{6}
\int_0^1
d F(t) \,Y(t) = c
\ ,
\ee
where $ F $ is an $ (n \times n) $-matrix
whose components are functions of bounded variation.

When the right-hand side of (\ref{4}) is
a continuous vector function $ g $,
it is well known that the  system
\be
\label{deterministic}
\left\{
\begin{array}{l}
{\ds
D Y(t) + A(t) Y(t) =
g(t)
\ ,\quad
t \in [0,1]
}
 \\ [2mm]
 \ds \int_0^1
 d F(t) \,Y(t) = c
\end{array}
\right.
\ee
admits a unique solution,
which belongs to $ C^{1}([0,1];\rr^n) $, if
and only if for some $s\in [0,1]$ (equivalently,
for every $s\in [0,1]$)
\[
\mbox{(H0)}
\hs{5cm}
\det \int_0^1 d F(t) \,\Phi^s(t) \neq 0
\ ,
\hs{5.5cm}
\]
where $ \Phi^s(t) $
denotes the fundamental matrix solution of
$ {\ds D Y(t) + A(t) Y(t) = 0} $, that is,
$\forall s \in [0,1]$,
\[
\left\{
\begin{array}{l}
{\ds \frac{d}{dt} \Phi^s(t) + A(t) \Phi^s(t) = 0
\ ,\quad
t \in [0,1] }
 \\ [2mm]
\Phi^s(s) = \mbox{I}
\ ,
\end{array}
\right.
\]
with $ \mbox{I} $ the identity matrix.
In turn, this is equivalent to say that the homogeneous problem
($g\equiv 0$, $c\equiv 0$) has only the trivial solution.
When hypothesis (H0) holds, the solution to (\ref{deterministic})
is given by
\[
Y(t)  =
J(t)^{-1} c +
\int_0^1 G(t,s) g(s) \,d s
\ ,
\]
where
\be
\label{8.2}
J(t) = \int_0^1 d F(u) \, \Phi^t(u)
\ee
and $ G(t,s) $ is the (matrix-valued) Green function associated to
$A$ and $F$. An explicit expression for this
function is the following
(see e.g. \cite{co1} or \cite{hon}):
\be
\label{8.3}
G(t,s) = J(t)^{-1} 
\Big[ \int_0^s d F(u) J(u)^{-1} - \1_{\{t \leq s\} } \mbox{I} \Big] J(s)
\ .
\ee
Under (H0), we define the solution to
(\ref{4})-(\ref{a-fbcY}) as the $n$-dimensional stochastic
process
\be
\label{8}
Y(t) =
J(t)^{-1} c + \int_0^1 G(t,s) \, dB(s)
\ ,
\ee
and the solution to (\ref{1})-(\ref{2}) as the
process $\{X(t)=Y_n(t),\ t\in [0,1]\}$.
The Green function (\ref{8.3}) has bounded variation,
so that the Wiener integrals in (\ref{8})
can be interpreted
pathwise by means of an integration by parts
$$
  \Big[\int_0^1 G(t,s)\,dB(s)\Big](\omega)=
  -\int_0^1 G(t,ds)\,B(s)(\omega)
$$
(we take into account here that $G(t,1)=0$, $\forall t$), and therefore
$Y$ can be defined everywhere. We shall assume throughout the paper
that the solution is interpreted in this pathwise sense.
Furthermore, it is not difficult to verify that the process $Y(t)$
so defined is continuous (hence $X(t)$ is a $C^{n-1}$ process)
and that, for each $ t\in [0,1] $,
the mapping $\omega\mapsto Y(\omega)$
from $\Omega$ into $C([0,1];\rr^n)$ is continuous with the usual
topologies.

Notice that, with the notation introduced in (\ref{6}),
the particular lateral condition (\ref{a-fbcY}) corresponds to
\be
\label{6.2}
dF =
\left[
\begin{array}{ccccc}
0      & \cdots & 0      & \sum\limits_{j=1}^{m} \alpha_{1j} \delta_{t_j}
\\
0      & \cdots & 0      & \sum\limits_{j=1}^{m} \alpha_{2j} \delta_{t_j}
\\
\vdots &        & \vdots & \vdots
\\
0      & \cdots & 0      & \sum\limits_{j=1}^{m} \alpha_{nj} \delta_{t_j}
\end{array}
\right]
\ ,
\ee
where $ \delta_t $ denotes the Dirac measure at $ t $, and
that $J_{ik}(t)=\sum_{j=1}^{m} \alpha_{ij} \Phi^t_{nk}(t_j)$.
Notice also that only the first column of $G(t,s)$ is relevant in
(\ref{8}).
\par
\bigskip
  Another natural definition of solution for the system
  (\ref{4})-(\ref{a-fbcY})
  arises if,
  for each $\omega$ fixed, we consider the object $\dot B(\omega)$ as
  the derivative of a continuous
  $\rr^n$-valued function defined on $[0,1]$, and therefore we
  regard (\ref{4}) as an equation between distributions.
  The vector function $Y=(Y_1,\dots,Y_n)$ will be a solution in
  the distributional sense if for any smooth vector  
  $\varphi=(\varphi_1,\dots,\varphi_n)$ vanishing
  in the complement of $\mathopen]0,1\mathclose[$, with $\int_0^1 \varphi=0$,
\be
\label{a-distsol}
  \int_0^1 \Big(Y(t)+\int_0^t A(a)Y(s)\,ds -B(t)\Big)\cdot \varphi(t)\,dt=0
\ee
  and (\ref{a-fbcY}) is satisfied.
  But (\ref{a-distsol}) amounts to say that there exists a constant $Y(0)$ such that
$$
  Y(t)-Y(0)+\int_0^t A(s)Y(s)\,ds=B(t)
  \ ,\quad
  t\in[0,1]
  \ ,
$$
  and a fortiori we find 
  that $Y$ must be a continuous function. It is easily seen that 
  both concepts of solution coincide.
\subsection{On the law of the solution}
In the present subsection we shall prove an
absolute continuity result 
for the law of the solution process $ \{Y(t),\ t \in [0,1] \} $.
Here we allow the boundary condition (2.4) to depend on all coordinates
of $Y$, since we will use this generality later on. 
\par
  Let $Y=(Y_1,\dots,Y_n)\colon\Omega\rightarrow C([0,1];\rr^n)$ be the 
  solution to the problem
\begin{equation}\label{a-lawproblem}
\left\{
\begin{array}{l}
{\ds
DY(t)+A(t)Y(t)=\dot B(t)
\ ,\quad
t \in [0,1]
}
 \\ [2mm]
 \Lambda[Y]=c
 \ ,
\end{array}
\right.
\end{equation}
  where $\Lambda$ is any linear operator on $C([0,1];\rr^n)$ with 
  finite support $\supp\Lambda=\{t_1,\dots,t_m\}$. (We are not 
  assuming here that $\Lambda$ involves only the coordinate function $Y_n$, 
  but we do assume that problem (\ref{a-lawproblem}) is well-posed.)
\par
  If $\{s_1,\dots,s_k\}\subset [0,1]$ is a set containing $\supp\Lambda$, 
  then $\Lambda$ can be regarded as a linear operator on the space of 
  functions $\big(\{s_1,\dots,s_k\}\rightarrow\rr^n\big)\cong \rr^{n\times k}$.
  We keep the same symbol $\Lambda$ for both interpretations. 
  Denote by $M$ the linear manifold in $\rr^{n\times k}$:
\begin{equation} \label{a-spaceH}
  M:=\Big\{
  x=\Big(
  \left(
  \begin{array}{c}
    x_{11}
    \\
    \vdots
    \\
    x_{1n}
  \end{array}
  \right) 
  ,\dots,
  \left(
  \begin{array}{c}
    x_{k1}
    \\
    \vdots
    \\
    x_{kn}
  \end{array}
  \right) 
  \Big)
  \in\rr^{n\times k}
  :\ 
  \Lambda[x]=c
  \Big\}
  \ .
\end{equation}
\begin{pr}
\label{abscon}
  With the notations above, the random vector 
  $(Y(s_1),\dots, Y(s_k))\colon\Omega\rightarrow \rr^{n\times k}$
  is absolutely continuous with respect to the Hausdorff measure
  in $M$.
\end{pr}

\no
{\em Proof}: 
  Taking into account that the vector $(Y(s_1),\dots,Y(s_k))$
  is Gaussian, it suffices to prove that any open ball
  in $M$ has a positive probability under the law of this vector.
\par
  Fix $x\in M$. Let us see first that there exists $\omega\in\Omega$
  such that the function 
  $Y(\omega)\colon [0,1]\rightarrow\rr^n$
  satisfies 
  $(Y(\omega)(s_1),\dots,Y(\omega)(s_k))=x$.
  Indeed, by simple interpolation,
  there obviously exists a $C^{\infty}$ function
  $y\colon [0,1]\rightarrow\rr^n$ such that
  $y_i\equiv y_{i+1}'$, 
  $i=1,\dots,n-1$, and $(y(s_1),\dots,y(s_k))=x$ 
  (therefore $\Lambda[y]=c$). 
  Defining
\begin{equation}\label{a-omegaofy}
  \omega(t)=y^1(t)-y^1(0)+
  \int_0^t (a_{n-1}(s)y^1(s)+\cdots+a_0(s)y^n(s))\,ds
  \ ,
\end{equation}
  we find that $y$ is the solution path
  $Y(\omega)$ of (\ref{a-lawproblem}).
  Any open ball $U(x)$ of $M$ 
  centred at a point $x\in M$ has therefore a non-empty 
  inverse image $Y^{-1}(B(x))\subset \Omega$. 
  Moreover, since the mapping $\omega\mapsto Y(\omega)$ is continuous,
  $Y^{-1}(U(x))$ is open. We get that $P\{Y\in U(x)\}>0$.
\qed
\begin{re}\label{a-disjdom}
  {\rm
  Proposition \ref{abscon} remains valid, 
  with a similar proof, if
  the domain where the problem is considered consists of two disjoint
  intervals, say $[0,a]$ and $[b,1]$, instead of a single one 
  (in that case $2n$ lateral 
  conditions are necessary for the problem to be well-posed). 
  The function $\omega$ can be defined as in (\ref{a-omegaofy})    
  for $t\in [0,a]$; as 
\[
  \omega(t)=y^1(t)-y^1(b)+
  \int_b^t (a_{n-1}(s)y^1(s)+\cdots+a_0(s)y^n(s))\,ds
\]
  for $t\in[b,1]$, and arbitrarily (continuous) on 
  $\mathopen]a,b\mathclose[$.
\qed
  }
\end{re}
\subsection{An approximation result}
\label{anapres}
We shall now state an easy approximation result (Proposition \ref{4zzz})
that we will need partially in the proof of Proposition \ref{noassumeA}. 
Consider the space $C^k:=C^k([0,1];\rr^{n\times n})$, 
with $k$ a fixed non-negative
integer or $\infty$, endowed with its natural topology. 
  Let ${\cal C}$ be the subset of $C^k$ comprising the matrix functions
  $A\colon [0,1]\rightarrow\rr^{n\times n}$ of the form (\ref{5}), 
  with the topology induced by $C^k$.
  Fix a linear operator
  $ \Lambda\colon C([0,1];\rr^n)\rightarrow \rr^n$ 
  of rank $n$ and with the form given by (\ref{6.2}),
  and
  consider the deterministic problems:
\be
\label{OR1}
\left\{
\begin{array}{l}
D Y(t) + A(t) Y(t) = 0
 \\ [2mm]
\Lambda [Y] = 0
\end{array}
\right.
\ee
with $ A \in C^k$. Let ${\cal D}\subset C^k$ the class of matrix functions $A$ such that
(\ref{OR1}) has only the trivial solution. Finally set $V:={\cal C}\cap{\cal D}$. 


\begin{lm}
\label{1zzz}
  $V$ is open and dense in ${\cal C}$.
\end{lm}

\no
{\em Proof}: 
For $A\in{\cal D}$, denote by $\Phi_A^0(t)$ the fundamental matrix solution of
the system 
$DY(t)+A(t)Y(t)=0$, with $Y(0)=\mbox{I}$. Consider the composition
of linear operators
\composiciodeaplicacions:
  nom: {}
  conjunt1: $\rr^n$
  aplicacio1: $\Gamma_A$
  conjunt2: $C([0,1];\rr^n)$
  aplicacio2: $\Lambda$
  conjunt3: $\rr^n$
  element1: $c$
  element2: $[t\mapsto\Phi_A^0(t)c]$
  element3: $\Lambda[t\mapsto\Phi_A^0(t)c]$
  fi
The mapping
\aplicacio:
  nom: {}
  conjunt1: $C^k$
  aplicacio: $\varphi$
  conjunt2: ${\cal L}(\rr^n;\rr^n)$
  element1: $A$
  element2: $\Lambda\circ\Gamma_A$
  fi
is continuous: Indeed,
$$
  \|{\Lambda\circ\Gamma_A-\Lambda\circ\Gamma_B}\|_{
  {\cal L}(\rr^n;\rr^n)}
  \le
  \|{\Lambda}_{{\cal L}(C,\rr^n)}\|\cdot
  \|{\Gamma_A-\Gamma_B}\|_{{\cal L}(\rr^n,C)}
  \ ,
$$
  and
$$
  \|{\Gamma_A-\Gamma_B}\|_{{\cal L}(\rr^n,C)}
  =
  \sup_{\|{c}\|=1}
  \|{\big(\Phi_A^0(t)-\Phi_B^0(t)\big)c}\|_{\infty}
  =\max_{i,j} \sup_t
  |\Phi_A^0(t)_{i,j}-\Phi_B^0(t)_{i,j}|
  \ .
$$
  The continuity follows from the uniform continuous dependence
  of the solution with respect to the data in a linear Cauchy problem.
\par
  Now we use the fact that the set $H$ of invertible operators on
  $\rr^n$ is open in ${\cal L}(\rr^n,\rr^n)$. We obtain that
  $\varphi^{-1}(H)$ is open in
  $C^k$. However
  $\varphi^{-1}(H)$ is the set of matrices $A$ such that
  $\det(\Lambda\circ\Gamma_A)\neq 0$, which coincides with
  $\cal D$ by definition. 
  This shows that $\cal D$ is open in $C^k$. 
\par
  Since $\cal C$ is a linear 
  manifold in $C^k$, we have that $V={\cal C}\cap{\cal D}$ is open 
  in ${\cal C}$. Note that the particular form of $\Lambda$ 
  does not play any role up to this point.
\par
  To prove the density, we start by checking that $V\neq\emptyset$.
  If $ A \in {\cal C}$ and $s\in[0,1]$, the corresponding
fundamental matrix $\Phi^s_A(\cdot)$ has the form
\[
  \Phi^s_A(t)=
\left[
\begin{array}{cccc}
D^{n-1}\phi_1 (t) & \cdots & D^{n-1}\phi_{n} (t)
\\
\vdots & & \vdots
\\
D\phi_1 (t) & \cdots & D\phi_{n} (t)
\\
\phi_1 (t) & \cdots & \phi_{n} (t)
\end{array}
\right]
\]
with
$$
  \Phi^s_A(s)=\mbox{I}
  \ ,
$$
for some $ C^{n+k}$ real functions $ \phi_1, \ldots, \phi_n $.
Conversely, any such matrix is the fundamental matrix solution
$\Phi^s$ of $DY(t)+A(t)Y(t)=0$ for some $A\in {\cal C}$.
We have
\[
\int_0^1 d F(t) \,\Phi^s(t) =
\left[
\begin{array}{cccc}
\alpha_{11} & \cdots & \alpha_{1m}
\\
\vdots & & \vdots
\\
\alpha_{n1} & \cdots & \alpha_{nm}
\end{array}
\right]
\left[
\begin{array}{cccc}
\phi_1 (t_1) & \cdots & \phi_{n} (t_1)
\\
\vdots & & \vdots
\\
\phi_1 (t_m) & \cdots & \phi_{n} (t_m)
\end{array}
\right]
\ .
\]
Take $s\not\in\{t_1,\dots,t_m\}$.
Since $ (\alpha_{ij}) $ has full rank and
$ n \le m $, we can obviously find numbers
$ \phi_i(t_j) $, $ 1\le i\le n$, $ 1\le j\le m $
such that this product is an invertible square matrix.
Then we take $ C^{n+k}$ functions $ \phi_1, \ldots, \phi_n $
interpolating these numbers and so that
$\big(D^{n-j}\phi_i(s)\big)_{i,j}=\mbox{I}$.
The corresponding $A$ will therefore
belong to $V$.

Given now $ A \in {\cal D} $, let us fix
$ A_0 \in V $.
For $ \lambda \in \rr $, define
\[
M_\lambda :=
\left[
\begin{array}{ccccc}
0 & 0 & \cdots & 0
\\
-1 & 0 & \cdots & 0
\\
\vdots& \ddots & \ddots & \vdots
\\
0 & \cdots & -1 & 0
\end{array}
\right]
+ \left[
\begin{array}{ccccc}
\lambda & 0 & \cdots & 0
\\
0 & 0 & \cdots & 0
\\
\vdots & \vdots & \ddots & \vdots
\\
0 & 0 & \cdots & 0
\end{array}
\right] A_0 + \left[
\begin{array}{ccccc}
1-\lambda & 0 & \cdots & 0
\\
0 & 0 & \cdots & 0
\\
\vdots & \vdots & \ddots & \vdots
\\
0 & 0 & \cdots & 0
\end{array}
\right] A
\ .
\]
Let us see that the function
$u\colon\lambda\mapsto \det\Lambda[\Phi^0_{M_{\lambda}}]$ 
is analytical: Indeed, $M_{\lambda}$ depends analytically
on $\lambda$, and so the fundamental solution $\Phi^0_{M_{\lambda}}$
is also analytic in $\lambda$.
Finally, analyticity is preserved by
the linear functional $\Lambda$ and the determinant.
Now assume $u\equiv 0$ in a neighbourhood of 0. This would imply
$u\equiv 0$ on the whole line. However
$$
  u(1)=\det\Lambda[\Phi_{A_0}^0]\neq 0\ .
$$
We conclude that there exists a sequence $\{\lambda_n\}_n$
converging to zero such that $M_{\lambda_n}\in V$.
  Since 
  $M_{\lambda_n}\to A$ as $n\to\infty$,
  the density is proved.
\par
  This proof borrows some ideas from Theorem 7.1 in Chow and Lasota
  \cite{ChowLaso}.
\qed
\par
\bigskip
\begin{pr}
\label{4zzz}
Let $ A^N(t) $ be a sequence of functions in $V$
converging to $ A(t) \in V$.
Let $ Y^N(t) $ and $ Y(t) $ be the corresponding unique solutions
to:
\[
\left\{
\begin{array}{l}
D Y^N(t) + A^N(t) Y^N(t) = \dot{W}(t)
 \\ [2mm]
\Lambda [Y] = c
\end{array}
\right.
\quad
\mbox{and}
\quad\quad
\left\{
\begin{array}{l}
D Y(t) + A(t) Y(t) = \dot{W}(t)
 \\ [2mm]
\Lambda [Y] = c
\ .
\end{array}
\right.
\]
Then $ Y^N(t) $ converges to $ Y(t) $ pointwise and in
$ L^p $, for all $ p \geq 1 $, uniformly in $ t $, that is:
\[
\lim_{N \rightarrow + \infty} \,
\sup_{0\le t\le 1} \big| Y^N(t) (\omega) - Y(t)(\omega) \big|
= 0
\ ,\quad \forall \omega \in \Omega \ ,
\]
\[
\lim_{N \rightarrow + \infty} \,
\sup_{0\le t\le 1} \|{Y^N(t) - Y(t)}\|_{L^p(\Omega)}
= 0
\ .
\]
\end{pr}

\no
{\em Proof}: 
  In the situation given, the fundamental solutions $\Phi^t(u)^N$
  converge to the fundamental solution $\Phi^t(u)$ uniformly in
  $t$ and $u$. From this fact one shows easily that $J^N$ and $(J^N)^{-1}$
  defined by (\ref{8.2}) converge uniformly to $J$ and $J^{-1}$, taking into account that the
  entries of $dF$ are finite measures. Hence, the Green functions
  $G^N(t,s)$ converge to $G(t,s)$ uniformly in $t$ and $s$ as well.
\par
  We have
$$
  \|{Y^N(t)-Y(t)}\|_{L^p}
  \le
  \big| J^N(t)^{-1} c-J(t)^{-1}c\big|
  +
  \Big\|{\int_0^1 \big(G^N(t,s)-G(t,s)\big)\,dB_s}\Big\|_{L^p}
  \ .
$$
  The first term tends to zero uniformly in $t$. For the second,
  note that
\begin{eqnarray*}
  & \ds
  \sup_t \mbox{E}\Big[\Big(\int_0^1
  \big(G^N_{i,j} (t,s)-G_{i,j}(t,s)\big)\,dW_s\Big)^2\Big]
  =
  \sup_t \int_0^1
  \big(G^N_{i,j} (t,s)-G_{i,j}(t,s)\big)^2\,ds
  \\ & \ds
  \le
  \int_0^1 \sup_t
  \big(G^N_{i,j} (t,s)-G_{i,j}(t,s)\big)^2\,ds
  \le
  \sup_t \sup_s
  \big(G^N_{i,j} (t,s)-G_{i,j}(t,s)\big)^2
  \ ,
\end{eqnarray*}
  which converges to zero. Since all random variables are Gaussian,
  the convergence to zero of the second moments (uniformly in $t$)
  implies the convergence to zero of all moments, also uniformly in
  $t$. We have proved the second statement of the Theorem.
\par
  We turn to the pointwise convergence:
  Since $G^N_{i,j}(t,s)-G_{i,j}(t,s)$ is a function of bounded
  variation which tends to zero uniformly in $t$ and $s$, the finite
  measures $G^N_{i,j}(t,ds)-G_{i,j}(t,ds)$ tend weakly to zero,
  uniformly in $t$, and we have
$$
  \Big|\int_0^1
  (G^N_{i,j} (t,s)\,dW_s
  -
  G_{i,j} (t,s)\,dW_s)\Big|
  =
  \Big|\int_0^1
  W_s\cdot \big(G^N_{i,j} (t,ds)
  -
  G_{i,j}(t,ds)\big)\Big|
  \to 0
  \ .
$$                                                    
\qed
\setcounter{equation}{0}
\section{A characterisation of conditional independence and convergence of
         conditional expectations}
\label{ccicce}
  In this section we state two facts of a general nature that will be our main 
  probabilistic tools in Section \ref{mplfbvp}. Lemma \ref{2g} is an abstract 
  result on the conditional independence of two random vectors when a function
  of them (of a special structure) is given; it was proved in \cite{afn} 
  (see also \cite{fn1}). Proposition \ref{gogknud}, 
  on the other hand, provides a sufficient condition for the $L^2$-convergence
  as $N\to\infty$ of a sequence of conditional expectations of the form 
  $
  \E [ F(U_1^N) | U_2^N ]
  .
  $
\par
  We will mention first 
  three auxiliary lemmas on the
  conditional independence of $ \sigma $-fields,
  whose proofs are not difficult.
  Recall that we write
$
{\cal F}_1 \ci\limits_{{\cal G}}
{\cal F}_2
$
to mean that the $ \sigma $-fields
$ {\cal F}_1 $ and $ {\cal F}_2 $ are
conditionally independent
given the $ \sigma $-field $ {\cal G} $.
\begin{lm}
\label{2d} 
Let ${\cal F}_1, {\cal F}_2, {\cal G}, {\cal F}'_1, {\cal F}'_2$
be $\sigma$-fields such that
${\cal F}'_1\subset{\cal F}_1\vee{\cal G}$ and
${\cal F}'_2\subset{\cal F}_2\vee{\cal G}$.
Then,
\[
{\cal F}_1 \ci_{{\cal G}}
{\cal F}_2
\quad
\Rightarrow
\quad
{\cal F}^\prime_1 \ci_{{\cal G}}
{\cal F}^\prime_2
\ .
\]
\qed
\end{lm}

\begin{lm}
\label{222d} 
Let ${\cal F}_1, {\cal F}_2, {\cal F}_3 $ and ${\cal G}$
be $\sigma$-fields such that
$
{\cal F}_1\vee{\cal F}_2 \ci\limits_{{\cal G}}
{\cal F}_3$
and
${\cal F}_1 \ci\limits_{{\cal G}}
{\cal F}_2
$.
Then, 
\[
{\cal F}_1\vee{\cal F}_3 \ci_{{\cal G}}
{\cal F}_2
\quad
\mbox{and}
\quad
{\cal F}_2\vee{\cal F}_3 \ci_{{\cal G}}
{\cal F}_1
\ .
\]
\qed
\end{lm}
\begin{lm}
\label{22d} 
Let ${\cal F}_1, {\cal F}_2$ and ${\cal G}$
be $\sigma$-fields such that
$
{\cal F}_1 \ci\limits_{{\cal G}}
{\cal F}_2
$
and ${\cal G} \subset {\cal F}_1$. 
Then, for any $\sigma$-field ${\cal H}$, with 
${\cal G} \subset {\cal H}\subset{\cal F}_1$,
\[
{\cal F}_1 \ci_{{\cal H}}
{\cal F}_2
\ .
\]
\qed
\end{lm}
Let
$( \Omega, {\cal F}, P ) $
be a probability space and
$ \cF_1 $ and $ \cF _2 $ two independent
sub-$ \sigma $-fields of $ \cF  $.
Consider two functions
$g_1 : \rr^d \times \Omega \rightarrow \rr^d $ and
$g_2 : \rr^d \times \Omega \rightarrow \rr^d $
such that
$ g_i $ is
$   {\cal B} (\rr^d) \otimes \cF_i $--measurable,
$  i = 1, 2 $.
Set
$ B(\varepsilon) := \{ x \in \rr^d , | x | < \varepsilon \} $,
and denote by $ \lambda $ the Lebesgue measure on $ \rr^d $.
Let us introduce the
following hypotheses:

\bdes

\item[(H1)]
There exists $\varepsilon_0>0$ such that
for almost all $\omega \in \Omega$, and for any
$   | \xi | < \varepsilon_0   $,
$   | \eta | < \varepsilon_0   $ the system
\[
\left\{
\begin{array}{l}
z_1 - g_1 (z_2, \omega)   = \xi
 \\ [2mm]
z_2 - g_2 (z_1, \omega)   = \eta
\end{array}
\right.
\]
has a unique solution
$ (z_1,z_2) \in \rr^{2 d} $.

\item[(H2)]
For every
$ z_1 \in \rr^d $ and $ z_2 \in \rr^d $,
the random vectors
$ g_1(z_2, \cdot) $ and $ g_2(z_1, \cdot) $
possess absolutely continuous
distributions and   the
function
\[
\delta (z_1,z_2) =  \sup_{0 < \varepsilon < \varepsilon_0}
\frac{1}{\lambda (B(\varepsilon))^2} \, P \{
\left|   z_1 - g_1(z_2)   \right| < \varepsilon ,\ 
\left|   z_2 - g_2(z_1)   \right| < \varepsilon   \}
\]
is locally integrable in $ \rr^{2 d} $, for some $\varepsilon_0>0$.

\item[(H3)]
For almost all
$   \omega \in \Omega$, the functions
$ z_2 \mapsto g_1 (z_2,\omega)$ and
$z_1 \mapsto g_2 (z_1,\omega)$
are continuously differentiable and
\[
\sup_{
 {{\sss | z_2 - g_2 (z_1, \omega) | < \varepsilon_0 }  \atop
 {\sss | z_1 - g_1 (z_2, \omega) | < \varepsilon_0 }}
 }
\big|
\det [ \mbox{I} - \nabla g_1(z_2,\omega)
\nabla g_2 (z_1,\omega) ]\big|^{-1} \in L^1(\Omega)
\]
for some $\varepsilon_0>0$,
where $\nabla g_i$ denotes the Jacobian matrix  of $g_i$
with respect to the first argument.
\edes

\no
Note that hypothesis (H1) implies the
existence of two random vectors
$Z_1$ and $Z_2$
determined by the system
\[
\left\{
\begin{array}{l}
Z_1 (\omega)   =   g_1 ( Z_2 (\omega), \omega)  \\ [2mm]
Z_2 (\omega)   =   g_2 ( Z_1 (\omega), \omega)
\end{array}
\right.
\]

\begin{lm}
\label{2g}
(Alabert, Ferrante, Nualart \cite{afn}).
Suppose the functions $   g_1  $ and $   g_2   $ satisfy
the above hypotheses (H1) to (H3).
Then the following statements are equivalent:

\bdes
\itemsep=-3pt
\item[(i)]   $\cF_1$ and $\cF_2$ are conditionally independent
given the random vectors $ Z_1 , Z_2 $.
\item[(ii)]
There exist two functions
$ F_i : \rr^{2 d} \times \Omega \rightarrow \rr $, $ i = 1, 2 $,
which are
$ \cB(\rr^{2 d}) \otimes \cF_i $--measurable,
such that
\[
\Big|
\det [ \mbox{\rm I} -
\nabla g_1(Z_2) \nabla g_2 (Z_1)
] \Big|     =
F_1 (Z_1,Z_2,\omega)
 F_2 (Z_1,Z_2,\omega)
\ ,\quad
\mbox{a.s.}
\]
\edes
\qed
\end{lm}
\par
Goggin \cite{gog} gives
a sufficient condition for the 
convergence in distribution of a sequence of conditional expectations
of the form $\E [ F(U_1^N) | U_2^N ]$.
We reproduce here a slightly simplified version.
Combining this result with Lemma \ref{knud}, due to Knudsen
\cite{kn1}, we can easily prove our Proposition \ref{gogknud}.
\begin{lm}
\label{gogi}
(Goggin \cite{gog}).
Let $ U_1^N$ and $U_2^N$ be two sequences of random vectors
on a
probability space $ (\Omega, {\cal F}, \mbox{P}) $,
such that
$ (U_1^N,U_2^N) \longrightarrow (U_1,U_2)$ as $N\to\infty$ in distribution.
Assume that:
\begin{enumerate}
\itemsep=-3pt           
\item
There exist probabilities $Q^N $ on $(\Omega,\cal F)$ such that
$ P \ll Q^N $ on $ \sigma\{U_1^N,U_2^N\} $
and $ U_1^N $ and $ U_2^N $ are
independent under $Q^N $.
Denote
$ \ell^N(U_1^N,U_2^N) := \frac{\textstyle dP}{\textstyle dQ^N}$\ .
\item
There exists a probability $Q$ on $(\Omega,\cal F)$ under which
$ U_1 $ and $ U_2 $ are independent.

\item
The $Q^N $--distribution
of $ (U_1^N,U_2^N,\ell^N(U_1^N,U_2^N)) $ converges
weakly to the $Q$--distribution
of $ (U_1,U_2,\ell(U_1,U_2)) $, where $\ell$ is such that $\E_Q [\ell(U_1,U_2)] = 1 $.
\end{enumerate}

\no
Then:
\begin{enumerate}
\itemsep=-3pt           
\item
$P \ll Q $ on $ \sigma\{U_1,U_2\} $ and
$ \frac{\textstyle dP}{\textstyle dQ} = \ell(U_1,U_2) $;

\item
For every bounded continuous function
$ F $,
\[
\E_P \Big[F(U_1^N)|U_2^N\Big]
\longrightarrow E_P \Big[F(U_1)|U_2\Big]
\mbox{\ in distribution.}
\]
\end{enumerate}
\qed
\end{lm}

\begin{lm}
\label{knud}
(Knudsen \cite{kn1}).
Let $U_1^N$ and $U_2^N$ be two sequences of random vectors on a probability
space $(\Omega,\cal F, P)$.
Assume that, as $N\to\infty$,
\begin{enumerate}
\itemsep=-3pt           
\item
$ U_2^N \stackrel{P}{\longrightarrow} U_2 $.
\item
$ U_1^N \stackrel{L^2}{\longrightarrow} U_1 $, with
$ U_1 \in L^p$, for some $p > 2 $.
\item
$ \| \E [ U_1 | U_2^N ] \|_{L^2}
\longrightarrow
\| \E [ U_1 | U_2 ] \|_{L^2} $.
\end{enumerate}

\no
Then
\[
\E \Big[ U_1^N | U_2^N \Big]
\stackrel{L^2}{\longrightarrow}
\E \Big[ U_1| U_2 \Big]
\]
as $N\to\infty$.
\qed
\end{lm}

\no
  Combining Lemma \ref{gogi} and \ref{knud}, we get the
  following proposition.
\begin{pr}
\label{gogknud}
Let $U_1^N$ and $U_2^N$ be two sequences of
random vectors. Assume that, as $N\to\infty$,
\begin{enumerate}
\itemsep=-3pt           
\item
$ U_1^N \stackrel{P}{\longrightarrow} U_1 $ and
$ U_2^N \stackrel{P}{\longrightarrow} U_2 $.
\item
Hypotheses 1,2 and 3 of Lemma \ref{gogi} hold true.
\end{enumerate}

\no
Then, for any bounded and continuous function $ F $,
\[
\E \Big[ F(U_1^N) | U_2^N \Big]
\stackrel{L^2}{\longrightarrow}
\E \Big[ F(U_1) | U_2 \Big]
\ .
\]
\end{pr}

\no
{\em Proof}: 
Applying Lemma \ref{gogi},
we have that, for every bounded and continuous
$ F $,
\be
\label{4321}
\E \Big[ F(U_1^N) | U_2^N \Big]
\longrightarrow
\E \Big[ F(U_1) | U_2 \Big]
\mbox{\ in distribution.}
\ee
  From (\ref{4321}) and the fact that
$ F(U_1^N), F(U_1) \in L^\infty $,
we obtain the convergence of the $ L^2 $ norms:
\begin{equation}
\label{convnorms}
\| \E [ F(U_1^N) | U_2^N ] \|_{L^2}
\longrightarrow
\| \E [ F(U_1) | U_2 ] \|_{L^2}
\ .
\end{equation}
  On the other hand, since $F$ is bounded and 
  $U_1^N \stackrel{P}{\longrightarrow} U_1$, we also have
\begin{equation}
\label{L2conv}
  F(U_1^N) \stackrel{L^2}{\longrightarrow} F(U_1)
\ .
\end{equation}
  Now, (\ref{convnorms}) and (\ref{L2conv}) imply that
\[
\| \E [ F(U_1) | U_2^N ] \|_{L^2}
\longrightarrow
\| \E [ F(U_1) | U_2 ] \|_{L^2}
\]
  and we get the conclusion applying Lemma \ref{knud}.
\qed
\setcounter{equation}{0}
\section{Markovian properties of linear functional boundary value problems}
\label{mplfbvp}
In the study of boundary value stochastic problems,
one of the main interests has been to seek 
conditions on the coefficients
for the solution process to
satisfy some suitably defined Markov-type property.
Intuition suggests that a relation $h(X(0),X(1))=0$
will possibly prevent the Markov process property from holding in general.
One might think
that nevertheless the Markov field property,
which is defined below, will be satisfied.
It is easy to see that any Markov process is a Markov field
(see Jamison \cite{jam} for the continuous case and 
Alabert and Marmolejo \cite{am1} for a simple proof in the general
case). 
The converse is not true.
For instance, the processes $X(t)=W(t)-\alpha W(1)$ are Markov fields;
they are not Markov processes, except for the cases $\alpha=0$ and $\alpha=1$.
\begin{de}
\label{3a}
A process
$ \{ X(t),\ t \in [0,1] \} $ is said to be
a {\em Markov field\/} if for any
$ 0 \leq a < b \leq 1 $, the
$ \sigma $-fields
$ \sigma \{ X(t), t \in [a,b] \} $ and
$ \sigma \{ X(t), t \in \mathopen]a, b\mathclose[^c\} $ are
conditionally independent given
$ \sigma \{ X(a), X(b) \} $.
\qed
\end{de}
However, even this weaker property holds only in special
cases. For instance, in \cite{afn} it was shown that
the solution to
\be
\left\{
\begin{array}{l}
\dot X(t)=b(X(t)) + \sigma(X(t)) \circ \dot{W} (t)
\ ,\quad
t\in [0,1]
 \\ [2mm]
X(0)=\psi(X(1))
\ ,
\end{array}
\right.
\ee
where the stochastic integral is understood in the Stratonovich sense,
is a Markov field if and only if 
$\ds b(x)=A\sigma(x)+B\sigma(x)
\int_c^x \frac{1}{\sigma(t)}\,dt$, for some constants $A,B,c$.
As a corollary, in case
$\sigma$ is a constant (additive noise),
$X$ is a Markov field if and only if
$b$ is an affine function.

Our aim is to study the linear--additive case
when the additional condition
takes into account the value of the
solution in some interior points of the time interval.
The following simple example illustrates
that the situation changes.
\begin{example}
{\rm
Consider the first order system
\[
\left\{
\begin{array}{l}
{\ds \dot X(t) = \dot W(t)
\ ,\quad
t \in [0,1]}
 \\ [2mm]
{X({\frac{1}{2}}) + X(1) = 0 \ .}
\end{array}
\right.
\]
The solution is the process
\[
\textstyle X(t)
 = - \frac{1}{2}\big(W(\frac{1}{2}) + W(1)\big) + W(t)
\ ,
\]
which is not a Markov field.
Indeed, for $ a = 0 $ and $ {b = \frac{2}{3} } $,
the random variables
$ {X(\frac{1}{2}) } $ and $ X(1) $ are not
conditionally independent given
$ \sigma \{ X(a), X(b) \} $. Nevertheless, $X$ is a Markov field
when restricted to $[0,\frac{1}{2}]$ or $[\frac{1}{2},1]$.
\qed
}
\end{example}
\par
\bigskip
In order to formulate precisely the conditional independence property
enjoyed by the system (\ref{4})-(\ref{a-fbcY}), we introduce first 
some more concepts and notation.
\par
Let $\Lambda_1,\dots,\Lambda_n$ be the real-valued components of
a boundary operator $\Lambda$ of the form (\ref{a-fbcY}) and  
denote their support by 
$ \supp \Lambda_i :=\{t_j \in [0,1]:\ 
\alpha_{ij} \neq 0\} $. 
\begin{de}
\label{defpreserves}
We will say that $\Lambda_i$ {\em preserves\/} the pair $(a,b)$
if either
$\supp \Lambda_i\subset\mathopen]a,b\mathclose[$ or
$\supp \Lambda_i\subset[a,b]^c$. If this is true for all $i$
(i.e.  
there are no boundary conditions involving simultaneously points
inside and outside $[a,b]$),
then we 
will also say that $\Lambda$ {\em preserves\/} $(a,b)$.
\qed
\end{de}
\par
\bigskip
We want to prove that the solution $Y$ to the system
(\ref{4})-(\ref{a-fbcY}) satisfies the following conditional
independence property (Theorem \ref{mainth}): 
If $\Lambda$ preserves $(a,b)$, then 
\[
\sigma \{ Y(t), t \in [a,b] \}
\ci_{\sigma \{ Y(a), Y(b) \}}
\sigma \{ Y(t), t \in \mathopen]a,b\mathclose[^{c} \}
\ .
\]
More generally, this conditional independence is also true 
when $\Lambda$ does not preserve $(a,b)$, provided
the conditioning $\sigma$-field is enlarged with the variables $Y_n(t)$, for 
$t$ in $[a,b]$ and in the support of all non-preserving boundary operators
$\Lambda_i$
(Theorem \ref{main3}).
\par
  Since the boundary conditions can be written in many different
  equivalent ways, and the sets $\supp \Lambda_i$ (hence the 
  property of preserving an interval) depend on the representation chosen,
  we need, before proceeding further, some sort of ``canonical'' definition  
  of the linear operator $\Lambda$.
Given $ \supp \Lambda=\{t_1,\dots,t_m\}$, $\Lambda$ can be regarded
as a linear mapping $\rr^m\rightarrow\rr^n$, that means, an $n\times m$ matrix
acting on the vector $(Y_n(t_1),\dots,Y_n(t_m))$ 
(see notations of Section \ref{lsdewfbc}).
\par
A {\em basis\/} $B$ for an $n\times m$ matrix $\Lambda$ is any 
$n\times n$ minor with full rank.
For notational simplicity, assume that $B$ consists of the
firsts $n$ columns of $\Lambda$. 
Denoting by $N$ the non-basic columns, we can write  $\Lambda=(B,N)$.
Defining $\widetilde{\Lambda}=(I,B^{-1}N)$, 
the system 
of equations 
$\Lambda x=c$ can be written 
in the equivalent form $\widetilde{\Lambda}x=B^{-1}c$. 
In this situation, we shall say that $\widetilde{\Lambda}$ is a
{\em basic\/} expression of $\Lambda$ relative to the
basis $B$.  
In the following lemma we prove that this
representation can be considered ``canonical'' for our purposes,
since any pair $(a,b)$ will or will not be
preserved by any basic equivalent form of $\Lambda$.
  In the sequel, we will always assume, without explicit mention, 
  that the boundary condition is written in this form.
\begin{lm}
\label{canon}
Let $\widetilde{\Lambda}$ and $\widetilde{\Lambda}'$
be two basic expressions of $\Lambda$,
and
fix $0 \leq a < b \leq 1 $.
Then, $\widetilde{\Lambda}$ preserves $(a,b)$ if and only if   
$\widetilde{\Lambda}'$ preserves $(a,b)$.   
\end{lm}

\no
{\em Proof}: 
Without any loss of generality we can assume that
$\widetilde{\Lambda} = (\mbox{I}, N)$,
where $\mbox{I}$ is the $n \times n$ identity matrix and
$$
N=
\left[
\begin{array}{ccc}
\alpha_{1, n+1}&\cdots&\alpha_{1, m}
\\
\vdots&\ddots&\vdots
\\
\alpha_{n, n+1}&\cdots&\alpha_{n, m}
\end{array}
\right]
\ .
$$
All basic expressions of the original
matrix $\Lambda$ can be obtained by repeated Gaussian pivoting   
on entries of non-basic columns; when pivoting on $\alpha_{ik}$, 
the column $i$ 
leaves the basis (the identity matrix) and is replaced by 
column $k$.
Therefore, it is sufficient to prove 
the lemma for
$\widetilde{\Lambda}$ and a basic expression $\widetilde{\Lambda}'$ obtained from 
$\widetilde{\Lambda}$ by one pivoting operation.
\par
Let us assume that $\alpha_{1,n+1}\neq 0$ and that 
the operator
$\widetilde{\Lambda}_i$ does not preserve the pair $(a,b)$.
We are going to find an operator $\widetilde{\Lambda}'_j$ which neither preserves
$(a,b)$.
Rows $1$ and $i$ before and after 
pivoting on $\alpha_{1,n+1}$ are the following:
\[
\begin{array}{ccccccccccccccc}
&&&&&&i&&&&&&&&
\\
\widetilde{\Lambda}_1 = &
[ & 1&
0 & \cdots & 0 & 0 &
0 & \cdots & 0 & \alpha_{1,n+1} &
\alpha_{1,n+2} & \cdots & \alpha_{1,m} & ]
 \\ [2mm]
\widetilde{\Lambda}_i = &
[ & 0 &
0 & \cdots & 0 & 1 &
0 & \cdots & 0 & \alpha_{i,n+1} &
\alpha_{i,n+2} & \cdots & \alpha_{i,m} & ]
 \\ [6mm]
\widetilde{\Lambda}'_1 = &
[ & \frac{1}{\alpha_{1,n+1}}&
0 & \cdots & 0 & 0 &
0 & \cdots & 0 & 1 &
\frac{\alpha_{1,n+2}}{\alpha_{1,n+1}} & \cdots & 
\frac{\alpha_{1,m}}{\alpha_{1,n+1}} & ]
 \\ [2mm]
\widetilde{\Lambda}'_i = &
[ & -\beta_i &
0 & \cdots & 0 & 1 &
0 & \cdots & 0 & 0 &
\gamma_{i,n+2} & \cdots & \gamma_{i,m} & ]
\end{array}
\]
where ${\beta_i = \frac{\alpha_{i,n+1}}{\alpha_{1,n+1}}}$
and ${\ds\gamma_{i,j} = 
\alpha_{i,j} - \alpha_{1,j} \beta_i}$. 
\par
If $\widetilde{\Lambda}_1$ does not 
preserve $(a,b)$, then 
the result is trivially true, with $j=1$, 
since
$\widetilde{\Lambda}_1$ and $\widetilde{\Lambda}'_1$
have their non-zero coefficients in the same columns.
By a similar reason,
if $\beta_i = 0$ we can take $j=i$.
\par
  Assume finally that $\supp \widetilde{\Lambda}_1 \subset \mathopen]a,b\mathclose[$  
  (so that $i\neq 1$),
  and that $\beta_i \neq 0$. If 
  $t_i\in[a,b]^c$, the result is proved, since 
  $\widetilde{\Lambda}'_i$ links $t_1$ and $t_i$, and 
  $t_1\in\mathopen]a,b\mathclose[$; take $j=i$.
  If $t_i\in \mathopen]a,b\mathclose[$, then there exists 
  $ k \in \{n+1, \dots, m\} $
  such that $t_k\in[a,b]^c$ with 
  $\alpha_{i k}\neq0$;
  since $ \alpha_{1 k}=0 $ 
  (because $\supp \widetilde{\Lambda}_1 \subset \mathopen]a,b\mathclose[$), 
  we have
  $ \gamma_{i k}=\alpha_{i k}\neq0$ and we can take again $j=i$.
\qed
\par
\bigskip
  Now we can formulate our main result:
\begin{thm}
\label{mainth}
Suppose the system  
$$
\left\{
\begin{array}{l}
{\ds
DY(t)+A(t)Y(t)=\dot B(t)
\ ,
\quad 
t \in [0,1]
\ ,
}
 \\ [2mm]
\Lambda[Y]=c
\end{array}
\right.
$$
satisfies {\em (H0)},
and let
$Y=\{ Y(t),\  t \in [0,1] \} $
be its unique solution. Then, 
\be
\label{mainds}
\sigma \{ Y(t),\ t \in [a,b] \}
\ci_{\sigma \{ Y(a), Y(b) \}}
\sigma \{ Y(t),\ t \in \mathopen]a,b\mathclose[^{c} \}
\ee
if and only if
the pair $(a,b)$ is preserved by $\Lambda$.
\qed
\end{thm}

\no
Our main tool for the proof of the `if' part in Theorem \ref{mainth}
will be Lemma \ref{2g}.
The idea is the following: We will split the $2n$-dimensional
random vector $ (Y(a),Y(b)) $ into two vectors $ Z^1 $ and $ Z^2 $
of suitable dimensions, in such a way that $Z^1$
be a function of
$ Z^2 $ and the increments
of the Wiener process $ W $ in $ [a,b] $,
and in turn $ Z^2 $ be a function of $ Z^1 $
and the increments of $ W $ in $ \mathopen]a,b\mathclose[^{c}$. These mappings   
will play the role of
$ g_1 $ and $ g_2 $ in the set of hypotheses (H1) to (H3).
The first will be defined through the solution to equation
$DY(t)+A(t)Y(t)=\dot B(t)$, with the components of $Z^2$ fixed to a constant;
the second will be defined similarly, fixing the components of $Z^1$ to
a constant. However, this means that we need to solve
our differential equation with several sets of constraints,
which are different from the original set, and therefore we cannot
ensure a priori that these problems are well-posed. Consequently,
the above functions $g_1$ and $g_2$ need not exist
in general.
\par
  To solve this technical difficulty, we will resort to a two-step
  procedure. First, we will assume that all functional boundary value
  problems that we need to solve are indeed well-posed. Then,
  hypotheses (H1) to (H3) can be checked, and Lemma \ref{2g}
  applies directly, yielding the desired result. This is the
  goal of Proposition \ref{assumeA}. Secondly, we will use
  the approximation result of Subsection \ref{anapres} to show that
  the matrix $A$ can be approximated by perturbed matrices
  $A^N$ for which all boundary problems involved are well-posed
  and whose solutions $Y^N$ converge to the solution $Y$ of the
  original problem. Then, the convergence of conditional expectations
  given in Proposition \ref{gogknud} will allow to carry
  the conditional independence properties of $Y^N$ to the limit.
  This second step is the contents of Proposition \ref{noassumeA}.
  The `only if' part of the theorem is shown in Proposition \ref{onlyifpart}.
\par
  Let us formulate precisely the assumption needed for the first
  step: Set
\be
\label{zorro}
\begin{array}{l}
\ell = \#
\{ i :\ \supp \Lambda_i\subset [0,a\mathclose[\}\ ,
 \\ [2mm]
p = \#
\{ i :\ \supp \Lambda_i\subset \mathopen ]b,1]\} \ ,
 \\ [2mm]
q = \#
\{ i :\ \supp \Lambda_i\subset \mathopen ]a,b\mathclose[\} \ .
\end{array}
\ee
  We can assume that the equalities $\Lambda_i[X]=c_i$ are ordered
  in the following way:
\be
\label{ordercon}
\begin{array}{l}
\supp\Lambda_i\subset [0,a\mathclose[
\ ,\quad i=1,\dots,\ell
\ ,
 \\ [2mm]
\supp\Lambda_i\subset \mathopen]a,b\mathclose[
\ ,\quad i=\ell+1,\dots,\ell+q
\ ,
 \\ [2mm]
\supp\Lambda_i\subset \mathopen]b,1]
\ ,\quad i=\ell+q+1,\dots,\ell+q+p
\ ,
\end{array}
\ee
and the remaining equations (those involving points
both in $[0,a\mathclose[$ and in $\mathopen]b,1]$), carry the labels
$i=\ell+q+p+1,\dots,n$.

Consider now $DY(t)+A(t)Y(t)=0$ with the following sets of lateral conditions
and the specified domain:
\be
\label{zorro1}
\left\{
\begin{array}{l}
Y_j(a)  =  0
\ ,\quad
j=1,\dots,n-\ell
 \\ [2mm]
\Lambda_i[Y]=0
\ ,\quad
i=1,\dots,\ell
\end{array}
\right.
\quad,
\quad \mbox{on} \ [0,a] \ ;
\ee

\be
\label{zorro2}
\left\{
\begin{array}{l}
Y_j (b)  =  0
\ ,\quad
j=1,\dots,n-\ell-q
 \\ [2mm]
Y_j (a) =  0
\ ,\quad
j=n-\ell+1,\dots,n
 \\ [2mm]
\Lambda_i[Y]=0
\ ,\quad
i=\ell+1,\dots,\ell+q
\end{array}
\right.
\quad,
\quad \mbox{on} \ [a,b] \ ;
\ee

\be
\label{zorro3}
\left\{
\begin{array}{l}
Y_j (b)  =  0
\ ,\quad
j=n-\ell-q+1,\dots,n
 \\ [2mm]
\Lambda_i[Y]=0
\ ,\quad
i=\ell+q+1,\dots,\ell+q+p
 \\ [2mm]
\Lambda_i[Y]=0
\ ,\quad
i=\ell+q+p+1,\dots,n
 \\ [2mm]
{\ds Y_n (t) = 0
\ ,\quad \forall t\in\bigcup_{i=\ell+q+p+1}^n
\big(\supp\Lambda_i\cap[0,a]\big) }
\end{array}
\right.
\quad,
\quad \mbox{on} \ [b,1] \ ,
\ee
(notice that the third and fourth lines result in
$ n - \ell -q-p$ equations involving only points
in $ [b,1] $).
\begin{de} \label{regular}
We will say that the pair $(a,b)$ is {\em regular\/} 
if
$DY(t)+A(t)Y(t)=0$ together with any of the sets of conditions
(\ref{zorro1}), (\ref{zorro2}) or (\ref{zorro3}) has only the trivial solution.
Otherwise $(a,b)$ will be called {\em singular\/}.
\qed
\end{de}

\begin{pr}
\label{assumeA}
Suppose the system 
\be
\label{Ysystem}
\left\{
\begin{array}{l}
{\ds
DY(t)+A(t)Y(t)=\dot B(t)
\ ,
\quad 
t \in [0,1]
\ ,
}
 \\ [2mm]
\Lambda[Y]=c
\end{array}
\right.
\ee
satisfies (H0),
and let
$Y=\{ Y(t),\  t \in [0,1] \} $
be its unique solution.
Let $(a,b)$ be a regular pair preserved by $\Lambda$.
Then (\ref{mainds}) holds true.
\end{pr}

\no
{\em Proof}: 
Let us define the $ \sigma $-fields
\[
\begin{array}{l}
{\cal F}^i_{a,b} = \sigma \{ W_t - W_a ,\ t \in [a,b] \}
 \\ [2mm]
{\cal F}^e_{a,b} = \sigma \{ W_t ,\ t \in [0,a] \} \vee
\sigma \{ W_1 - W_t ,\ t \in [b,1] \}
\end{array}
\]
for $ 0 \leq a < b \leq 1 $.
Notice that $ {\cal F}^i_{a,b} $ and $ {\cal F}^e_{a,b} $
are independent.
\par
\medskip
We shall divide the proof
into several steps. In Step 1 we reduce the proof to that of
the conditional independence of two independent $\sigma$-fields.
In Step 2 it is shown that there exist the two functions
$g_1$ and $g_2$ needed to apply Lemma \ref{2g}. The hypotheses
of this lemma are checked in Steps 3, 4 and 5. In Step 6 we
finally conclude the result.

\medskip
\no
{\bf Step 1}
{\em Denote $ {\cal G}_{a,b} =
\sigma \{ Y(a), Y(b) \} $.
If
\be
\label{fifeg}
{\cal F}^i_{a,b}
\ci_{{\cal G}_{a,b}}
{\cal F}^e_{a,b} \ ,
\ee
then (\ref{mainds}) holds.}

\no
{\em Proof of Step 1}: 
It is immediate to prove that
$ \sigma \{ Y(t),\ t \in [a,b] \}
\subset {\cal G}_{a,b} \vee {\cal F}^i_{a,b} $, and
$ \sigma \{ Y(t),\ t \in \ \mathopen]a,b\mathclose[^{c} \}
\subset {\cal G}_{a,b} \vee {\cal F}^e_{a,b} $.
We apply then Lemma \ref{2d}.
\par
\medskip
\no
{\bf Step 2}
{\em Let $\ell$, $ p $ and $q$ be as in
$ (\ref{zorro}) $.
We will denote by
$ \widetilde{Y} $ the solution to (\ref{Ysystem}), to
distinguish the actual solution from $Y$ regarded as an unknown
of the system.
Define
\[
\begin{array}{l}
{\ds
Z^1 := 
\left( \widetilde{Y}_1(a), \ldots, \widetilde{Y}_{n-\ell}(a),
\widetilde{Y}_{n-\ell-q+1}(b), \ldots, \widetilde{Y}_n(b) \right)
\in\rr^{n+q}\ ,}
 \\ [2mm]
{\ds
Z^2 := 
\left( \widetilde{Y}_1(b), \ldots, \widetilde{Y}_{n-\ell-q}(b),
\widetilde{Y}_{n-\ell+1}(a), \ldots, \widetilde{Y}_n(a) \right)
\in\rr^{n-q}\ .}
\end{array}
\]
Then, there exist two functions
\[
\begin{array}{l}
{\ds
g_1\colon\rr^{n-q}\times\Omega\rightarrow\rr^{n+q}
\ ,}
 \\ [2mm]
{\ds
g_2\colon\rr^{n+q}\times\Omega\rightarrow\rr^{n-q}
\ ,}
\end{array}
\]
measurable with respect to
$ {\cal B} (\rr^{n-q}) \otimes {\cal F}^i_{a,b} $
and
$ {\cal B} (\rr^{n+q}) \otimes {\cal F}^e_{a,b} $
respectively, and
such that
\be
\label{26}
Z^1 = g_1 (Z^2, \omega)
\hs{0.5cm}
\mbox{and}
\hs{0.5cm}
Z^2 = g_2 (Z^1, \omega)
\ .
\ee}

\no
{\em Proof of Step 2}: 
Consider the lateral conditions
\be
\label{27}
\left\{
\begin{array}{l}
Y_j (b)  =  Z^2_j
\ ,\quad
j=1,\dots,n-\ell-q
 \\ [2mm]
Y_j (a) =  Z^2_{j-q}
\ ,\quad
j=n-\ell+1,\dots,n
 \\ [2mm]
\Lambda_i[Y]=
c_i
\ ,\quad
i=\ell+1,\dots,\ell+q
\end{array}
\right.
\ee
on $[a,b]$.
The process $ \widetilde{Y} $
trivially satisfies $DY(t)+A(t)Y(t)=\dot B(t)$ and these conditions on $[a,b]$;
however the solution to this problem is also unique.
Therefore, taking into account that $\Lambda_i[\tilde Y]=c_i$
are constants, the vector $Z^1$ is determined by $ Z^2 $
and the increments of the Wiener process in $ [a,b] $.
Moreover, the function $ g_1 (z_2,\omega) $ so defined
has a sense for every $ z_2 \in \rr^{n-q} $,
because we have that the solution to (\ref{4})-(\ref{27})
is unique, and
this fact does not depend on the particular right-hand sides.

We want to prove analogously the existence of the function
$ g_2 $. Consider first $DY(t)+A(t)Y(t)=\dot B(t)$ on $[0,a]$ with conditions
\be
\label{28}
\left\{
\begin{array}{l}
Y_j (a) = Z^1_j
\ ,\quad
j=1,\dots,n-\ell
 \\ [2mm]
\Lambda_i[Y]=c_i
\ ,\quad
i = 1, \ldots, \ell
\ .
\end{array}
\right.
\ee
The restriction to $ [0,a] $ of the solution
$ \widetilde{Y} $ to (\ref{Ysystem}) solves
also the differential system with conditions (\ref{28}), and is its unique solution.
Consider now
$DY(t)+A(t)Y(t)=\dot B(t)$ on $ [b,1] $, with
\be
\label{29}
\left\{
\begin{array}{l}
Y_j (b) = Z^1_{j+q}
\ ,\quad
j=n-\ell-q+1,\dots,n
 \\ [2mm]
\Lambda_i[Y]=c_i
\ ,\quad
i = \ell+q+1, \ldots, \ell+q+p
\end{array}
\right.
\ee
and the $ n-\ell-q-p $ equations on $ [b,1] $
that result from
\be
\label{291}
\left\{
\begin{array}{l}
\Lambda_i[Y]=c_i
\ ,\quad i = \ell+q+p+1, \ldots, n
 \\ [2mm]
{\ds Y_n (t) = \tilde{Y}_n (t)
\ ,\quad \forall t\in\bigcup_{i=\ell+q+p+1}^n
\big(\supp\Lambda_i\cap[0,a]\big)}
\ .
\end{array}
\right.
\ee
Again,
$ \widetilde{Y} $ restricted to $ [b,1] $
is its unique solution with conditions (\ref{29})-(\ref{291}).
The values $\tilde Y_n(t)$ appearing here are found in (\ref{29})
as a function of $Z^1_{j+q}$, $j=n-\ell-q+1,\dots,n$, and the Wiener
process on $[0,a]$. Therefore, the whole vector $Z^2$ is determined
by $Z^1$ and the increments of $W$ in $\mathopen]a,b\mathclose[^c$.
As before, the function $ g_2(z_1,\omega) $
so defined has a sense
for all $ z_1 \in \rr^{n+q} $.
\par
\medskip
\no
{\bf Step 3}
{\em
The functions $ g_1 $ and $ g_2 $ found in Step 2 satisfy (H1).}

\no
{\em Proof of Step 3}: 
The solution to a linear differential equation depends linearly
on the lateral data $c$ (see (\ref{8})).
Therefore, for each $ \omega $ fixed, system
(\ref{26}) is linear and it is enough to check that it has a
unique solution for $ \xi = \eta = 0 $.

Now, gathering together the lateral conditions (\ref{28}), (\ref{29}),
(\ref{291}), we obtain the original lateral conditions, so that system 
(\ref{26}) is equivalent to (\ref{Ysystem}) and therefore the
solution exists and is unique. 
\par
\medskip
\no
{\bf Step 4}
{\em $ g_1 $ and $ g_2 $ satisfy (H2).}

\no
{\em Proof of Step 4}: 
  The boundary value problem that defines $g_1$ consists of 
  equation $DY(t)+A(t)Y(t)=\dot B(t)$, together with the conditions (\ref{27}). 
  The resulting vector 
$$
  (Y_1(a),\dots,Y_{n-\ell}(a),Y_{n-\ell-q+1}(b),Y_n(b))
$$ 
  is 
  absolutely continuous on $\rr^{n+q}$, by Proposition \ref{abscon}. 
  The proof for $g_2$ is analogous, using Remark \ref{a-disjdom}. 

Finally, the random vectors $z_1-g_1(z_2,\omega)$ and $z_2-g_2(z_1,\omega)$
are independent and have the form $z_1-M^1z_2+U^1(\omega)$ and
$z_2-M^2z_1+U^2(\omega)$ respectively, for some constant matrices $M^1$ and $M^2$
and some Gaussian absolutely continuous vectors $U^1$ and $U^2$.
We deduce that the $\rr^{2n}$-valued random vector
$(z_1-g_1(z_2,\omega),z_2-g_2(z_1,\omega))$ has a density which is uniformly
bounded in $z_1$ and $z_2$. It follows at once that the function
$\delta$ in (H2) is bounded.
\par
\medskip
\no
{\bf Step 5}
{\em $g_1$ and $g_2$ satisfy (H3). Specifically,
$
\det [ \mbox{\rm I} -
\nabla g_1(z_2,\omega) \nabla g_2 (z_1,\omega)
]
$
is a constant different from zero.}

\no
{\em Proof of Step 5}: 
$ g_1 $ and $ g_2 $ are affine functions
of the first argument, with a
non-random linear coefficient (see (\ref{8})).
Therefore,
$ \nabla g_1(z_2,\omega) $ and $ \nabla g_2 (z_1,\omega) $
are constant matrices of dimensions $(n+q)\times(n-q)$ and
$(n-q)\times(n+q)$ respectively, which we denote simply
$ \nabla g_1 $ and $ \nabla g_2 $.
We know that the linear system
\[
\left\{
\begin{array}{l}
z_1   =   g_1 ( z_2, \omega)  \\ [2mm]
z_2   =   g_2 ( z_1, \omega)
\end{array}
\right.
\]
admits a unique solution.
This is equivalent to say
\[
\det [ \mbox{I} - \nabla g_1 \nabla g_2 ] =
\det 
\left(
\begin{array}{cc}
\mbox{I} & - \nabla g_1
\\
- \nabla g_2 & \mbox{I}
\end{array}
\right)
\neq 0
\ .
\]
\par
\medskip
\no
{\bf Step 6}
{\em
Relation (\ref{mainds}) holds true.}

\no
{\em Proof of Step 6}: 
We can apply Lemma \ref{2g} and the factorization in (ii)
trivially holds.
We deduce the relation (\ref{fifeg}) and,
by Step 1, that the process $Y$
satisfies the desired property,
for $(a,b)$ regular.
\qed
\par                 
\bigskip
  Let us now extend Proposition \ref{assumeA} to singular
  pairs $(a,b)$, using an approximation argument. We denote by 
  Eq($A,\Lambda$) our functional boundary value 
  problem relative to the matrix function $A$
  and the boundary operator $\Lambda$. The boundary data  
  $c$ will be fixed throughout. Let us call 
  $\Lambda_1$, $\Lambda_2$, $\Lambda_3$ the operators associated 
  to the lateral conditions 
  given by (\ref{zorro1}), (\ref{zorro2}) and (\ref{zorro3}), respectively. 
\begin{pr}
\label{noassumeA}
  Proposition \ref{assumeA} holds also for singular pairs $(a,b)$.
\end{pr}
\no
{\em Proof}: 
  Our initial hypothesis (H0) states that the original
  problem has one and only one solution, that is,
  $A(t)\in V_{\Lambda}$, where $V_{\Lambda}$ stands for the
  set $V$  defined in Subsection \ref{anapres}, relative to the boundary 
  operator $\Lambda$.
\par
  Fix $0\le a\le b$. 
  If the pair $(a,b)$ is singular, 
  then at least one of the
  problems Eq($A,\Lambda_1$), Eq($A,\Lambda_2$), or Eq($A,\Lambda_3$)
  is not well-posed. We know from Lemma \ref{1zzz}  
  that $V_{\Lambda}$
  is open and dense in the space ${\cal C}$ of matrices of the form (\ref{5}), 
  for any $\Lambda$. Therefore the set
  $V:=V_{\Lambda}\cap V_{\Lambda_1}\cap V_{\Lambda_2}\cap V_{\Lambda_3}$
  is also open and dense in ${\cal C}$.
\par
  Let $\{A^N(t),\ N\in\nn\}$ be a sequence of elements of $V$ converging to
  $A(t)\in{\cal C}$. From Proposition \ref{4zzz}, the corresponding solutions 
  $Y^N(t)$ converge to $Y(t)$
  in $L^p$.
\par
  Fix $s\in[a,b]^c$, and $r_1,\dots,r_k\in (a,b)$. 
  Consider the space $M$ defined in (\ref{a-spaceH}), based on the 
  coordinates $\{s,r_1,\dots,r_k\}\cup \supp\Lambda$. Let $M'$ be 
  the projection of $M$ onto the coordinates $r_1,\dots,r_k$. 
  Assume that $Y^N(s)$ are non-degenerate (if they are, $Y(s)$ will also be 
  a constant, and there is nothing to prove). Using Proposition
  \ref{abscon}, the vector 
  $U^N:=(Y^N(s),Y^N(r_1),\dots,Y^N(r_k))$ is a Gaussian vector with some 
  density $f^N$ 
  with respect to the Hausdorff measure on $\rr\times M'$. 
  Then there clearly exists an equivalent Gaussian probability 
  with density $f_0$  
  on $\rr\times M'$ whose first coordinate is incorrelated with 
  the remaining ones. 
  Define the probability
  $Q^N$ on $\Omega$ by $dP=\ell^NdQ^N$, with 
  $\ell^N=\frac{f^N}{f_0}(U^N)$. 
  Then hypothesis 1 of Lemma \ref{gogi} is clearly
  satisfied with $U_1^N=Y^N(s)$, $U_2^N=(Y^N(r_1),\dots,Y^N(r_k))$. 
  Analogously, one can define a probability $Q$
  by $dP=\ell dQ$, where 
  $\ell=\frac{f}{f_0}(U)$ and $f$ is the density of  
  $U:=(Y(s),Y(r_1),\dots,Y(r_k))$, satisfying hypothesis 2 of Lemma 
  \ref{gogi} with $U_1=Y(s)$ and 
  $U_2=(Y(r_1),\dots,Y(r_k))$. 
\par
  We prove now that hypothesis 3 also holds true: 
  we want to see that if $h$ is a bounded and continuous function, 
  then 
\begin{equation}\label{a-lawconv}
  \lim_{N\to\infty} \E_{Q^N}\big[h\big(U^N,\ell^N(U^N)\big)\big]
  =
  \E_{Q}\big[h\big(U,\ell(U)\big)\big]
  \ .
\end{equation}
  Notice that 
$$
  \E_{Q^N}\big[h\big(U^N,\ell^N(U^N)\big)\big]
  =
  \int_M h\big(u,\ell^N(u)\big)f_0(u)\,{\cal H}(du)
  \ ,
$$
  where ${\cal H}$ denotes the Hausdorff measure, 
  and analogously
$$
  \E_{Q}\big[h\big(U,\ell(U)\big)\big]
  =
  \int_M h\big(u,\ell(u)\big)f_0(u)\,{\cal H}(du)
  \ .
$$
  Convergence (\ref{a-lawconv}) is then easily derived through the 
  dominated convergence theorem.
  We can therefore apply Proposition \ref{gogknud} to obtain that
\begin{eqnarray*}
  & \ds
  L^2\mbox{--}\lim_N
  \E[F(Y^N(s))\ |\ Y^N(a), Y^N(r_1),
  \dots, Y^N(r_n), Y^N(b)]
  \\ & \ds
  =
  \E[F(Y(s))\ |\ Y(a), Y(r_1),
  \dots, Y(r_n), Y(b)]
  \ .
\end{eqnarray*}
  Analogously, one obtains the limit
$$
  L^2\mbox{--}\lim_N
  \E[F(Y^N(s))\ |\ Y^N(a),
  Y^N(b)]
  =
  \E[F(Y(s))\ |\ Y(a),
  Y(b)]
  \ .
$$
  We conclude that the conditional independence property
  can be carried to the limit and this finishes the proof of the present
  proposition and consequently of the necessity in Theorem \ref{mainth}.
  \qed 

\bigskip
  The `only if' part of Theorem \ref{mainth} is far easier to prove:
\begin{pr}
\label{onlyifpart}
  If $\Lambda$ does not preserve $(a,b)$, then relation (\ref{mainds}) is false.
\end{pr}
\no
{\em Proof:}    
Let us assume, to keep notations simple,
that there is only one boundary operator $\Lambda_1$,
which corresponds to a first order equation (the general
case can be stated similarly).
Assume that $\Lambda_1$ does not preserve $(a,b)$ and that (\ref{mainds}) 
holds true.
The corresponding boundary condition can thus be written as
\be
\label{boom1}
\sum_{t_j\in[a,b]^{c}} \alpha_j Y(t_j) + 
\sum_{t_j\in\mathopen]a,b\mathclose[} \alpha_j Y(t_j)=c
\ ,
\ee
where none of the summations is void.

\no
Notice first that
\be
\label{boom2}
\Big(\sum_{t_j\in\mathopen]a,b\mathclose[} 
\alpha_j Y(t_j), Y(a), Y(b)\Big)
\ee
is an absolutely continuous random vector
in $\rr^3$.
This follows easily from Proposition \ref{abscon}.

\no
Now, for any bounded and measurable function
$\Psi\colon\rr\rightarrow\rr$,
we have, by the conditional independence hypothesis and
relation (\ref{boom1}), that
$$
\begin{array}{l}
\mbox{E}[\Psi(\sum_{t_j\in\mathopen]a,b\mathclose[} 
\alpha_j Y(t_j))
| Y(a),Y(b)]=
\\ \\
=\mbox{E}[\Psi(\sum_{t_j\in\mathopen]a,b\mathclose[} 
\alpha_j Y(t_j))
| \{Y(t),\ t\in\mathopen]a,b\mathclose[^{c}\}]=
\Psi(\sum_{t_j\in\mathopen]a,b\mathclose[} \alpha_j Y(t_j))
\ .
\end{array}
$$
Taking $\Psi=\1_{[-M,M]}$ for some $M>0$,
we get that,
on $\{|\sum_{t_j\in\mathopen]a,b\mathclose[} 
\alpha_j Y(t_j)|\le M\}$,
which is a set of positive probability,
$ \sum_{t_j\in\mathopen]a,b\mathclose[} \alpha_j Y(t_j)$ is a 
measurable function of $(Y(a),Y(b))$.
In particular, this contradicts the absolute continuity of
(\ref{boom2}).
\qed
\par
\bigskip
  Next theorem generalises Theorem \ref{mainth} by allowing the existence of
  non-preserving boundary operators, at the price of enlarging the conditioning
  $\sigma$-field.
  The result can hardly be called a Markovian type property; nevertheless, 
  it seems interesting in itself, and gives rise to the conjecture 
  contained in Remark 
  \ref{conjecture} below. 
\begin{thm}
\label{main3}
Suppose the system 
$$
\left\{
\begin{array}{l}
{\ds
DY(t)+A(t)Y(t)=\dot B(t)
\ ,
\quad 
t \in [0,1]
\ ,
}
 \\ [2mm]
\Lambda[Y]=c
\end{array}
\right.
$$
satisfies (H0),
and let
$Y=\{ Y(t),\  t \in [0,1] \} $
be its unique solution. Fix $0\le a<b\le 1$ such that $a,b\notin \supp \Lambda$, and
let ${\cal G}$ be the $\sigma$-field generated by $Y(a)$, $Y(b)$ and all 
variables $Y_n(t)$, 
with
$t$ in $\supp\Lambda_i\cap \mathopen]a,b\mathclose[$ for some $\Lambda_i$ not preserving
$(a,b)$.
Then,
\be
\label{main3ds}
\sigma \{ Y(t), t \in [a,b] \}
\ci_{{\cal G}}
\sigma \{ Y(t), t \in \mathopen]a,b\mathclose[^{c} \}
\ .
\ee
  The same holds true replacing $\mathopen]a,b\mathclose[$ by $[a,b]^c$ in 
  the definition of ${\cal G}$.
\end{thm}
\no
{\em Proof:}   
  The result can be proved as Theorem \ref{mainth} with some modifications. We will only
  give a sketch of the necessary changes in the simple case where there is only
  one non-preserving boundary condition, which links one only point $t^*$ inside $[a,b]$
  with one or more points in $[0,a]$. Specifically, fix $a,b$ and assume
  that $\ell$, $p$ and $q$ are as in (\ref{zorro}) and (\ref{ordercon}), 
  the boundary conditions with support in $[0,a]\cup[b,1]$ carry the labels 
  $i=\ell+q+p+1,\dots,n-1$, and that
  $\Lambda_n$ is the non-preserving condition. 
\par
  In Step 2, define
\[
\begin{array}{l}
{\ds
Z^1 := 
\left( \widetilde{Y}_1(a), \ldots, \widetilde{Y}_{n-\ell-1}(a),
\widetilde{Y}_{n-\ell-q}(b), \ldots, \widetilde{Y}_n(b), \widetilde{Y}_n(t^*) \right)
\in\rr^{n+q+1}\ ,}
 \\ [2mm]
{\ds
Z^2 := 
\left( \widetilde{Y}_1(b), \ldots, \widetilde{Y}_{n-\ell-q-1}(b),
\widetilde{Y}_{n-\ell}(a), \ldots, \widetilde{Y}_n(a) \right)
\in\rr^{n-q}\ .}
\end{array}
\]
Consider the lateral conditions
\[
\left\{
\begin{array}{l}
Y_j (b)  =  Z^2_j
\ ,\quad
j=1,\dots,n-\ell-q-1
 \\ [2mm]
Y_j (a) =  Z^2_{j-q}
\ ,\quad
j=n-\ell,\dots,n
 \\ [2mm]
\Lambda_i[Y]=c_i
\ ,\quad
i=\ell+1,\dots,\ell+q
\end{array}
\right.
\]
on $[a,b]$. This system, as before, defines the function $g_1$. 
To define function $g_2$, we consider first 
$DY(t)+A(t)Y(t)=\dot B(t)$ on $[0,a]$ with conditions
\[
\left\{
\begin{array}{l}
Y_j (a) = Z^1_j
\ ,\quad
j=1,\dots,n-\ell-1
 \\ [2mm]
\Lambda_i[Y]=c_i
\ ,\quad
i = 1, \ldots, \ell
 \\ [2mm]
\Lambda_n[Y]=c_n
 \\ [2mm]
Y_n(t^*)=\widetilde{Y}_n(t^*)
\ ,
\end{array}
\right.
\]
and secondly
the system on $ [b,1] $, with
\[
\left\{
\begin{array}{l}
Y_j (b) = Z^1_{j+q}
\ ,\quad
j=n-\ell-q,\dots,n
 \\ [2mm]
\Lambda_i[Y]=c_i
\ ,\quad
i = \ell+q+1, \ldots, \ell+q+p
\end{array}
\right.
\]
and the $ n-1-\ell-q-p $ equations on $ [b,1] $
that result from
\[
\left\{
\begin{array}{l}
\Lambda_i[Y]=c_i
\ ,\quad i = \ell+q+p+1, \ldots, n-1
 \\ [2mm]
{\ds Y_n (t) = \tilde{Y}_n (t)
\ ,\quad \forall t\in\bigcup_{i=\ell+q+p+1}^{n-1}
\big(\supp\Lambda_i\cap[0,a]\big)}
\ .
\end{array}
\right.
\]
\par
  The claims of Step 4 are also easy to verify using Proposition 
  \ref{abscon}, taking into account that we assume $\Lambda$ is in basic 
  form, which implies that $\tilde Y_n(t^*)$ cannot be a constant. 
\qed
\par
\bigskip
  Another application of the ideas in the proof of Theorem \ref{mainth} provides
  the following ``Markov process'' property.
\begin{thm}
\label{main2}
  Fix $0\le a\le 1$ such that $a\notin \supp \Lambda$; the process $Y(t)$ satisfies
$$
\sigma \{ Y(t),\ t \in [0,a] \}
\ci_{\sigma \{ Y(a)\}}
\sigma \{ Y(t),\ t \in [a,1] \}
$$
  if and only if for all $i$, either
  $\supp\Lambda_i\subset[0,a\mathclose[$, or $\supp\Lambda_i\subset\mathopen]a,1]$.
\end{thm}
\no
{\em Proof}:   
  One can use 
  the same machinery as in the proof of Theorem \ref{mainth}:   
  The `only if' part can be proved within the same lines as Proposition \ref{onlyifpart},
  whereas 
  for the other implication,
  if  
  $\supp\Lambda_i\subset [0,a\mathclose[$,
  $i=1,\dots,\ell$, 
  and 
  $\supp\Lambda_i\subset \mathopen]a,1]$,
  $i=\ell+1,\dots,n$, then
  one can take 
\[
\begin{array}{l}
{\ds
Z^1 := 
( \widetilde{Y}_1(a), \ldots, \widetilde{Y}_{n-\ell}(a))
\in\rr^{n-\ell}\ ,}
 \\ [2mm]
{\ds
Z^2  := 
(\widetilde{Y}_{n-\ell+1}(a), \ldots, \widetilde{Y}_n(a))
\in\rr^{\ell}\ ,}
\end{array}
\]
  and define $g_1$ as a function of $Z^2$ and the increments of the Wiener process
  in $[0,a]$, and $g_2$ as a function of $Z^1$ and the increments of the Wiener
  process in $[a,1]$.
 \qed
\begin{re}
\label{conjecture}
{\rm
  Property
  (\ref{main3ds}) implies, using Lemma \ref{22d}, that  
\[
\sigma \{ Y(t), t \in [a,b] \}
\ci_{{\cal H}}
\sigma \{ Y(t), t \in \mathopen]a,b\mathclose[^{c} \}
\ ,
\]
  where ${\cal H}=\sigma\{Y(a),Y(b);\ Y(t),\,t\in\supp\Lambda\cap [a,b]\}$.
  We conjecture that this property holds true for a linear functional boundary
  operator $\Lambda$ supported on any subset of $[0,1]$, 
  and that it is false in general
  if ${\cal H}$ is replaced by a smaller $\sigma$-field.
  We will show a simple example illustrating the conjecture.
  Unfortunately,
  the technique we have employed here does
  not allow us to prove it. 
\par
  Consider the process $\ds X(t):=-\int_0^1 W(u)\,du+W(t)$, solution of the first order 
  problem
\[
\left\{
\begin{array}{l}
{\ds \dot X(t) = \dot W(t)
\ ,\quad
t \in [0,1]}
 \\ [2mm]
{\ds\int_0^1 X(u)\,du = 0 \ ,}
\end{array}
\right.
\]
  in which the support of the boundary operator is the whole interval $[0,1]$.
  Fix $a\in\mathopen]0,1\mathclose[$, and set 
\[
  {\cal G}:=
  \sigma \{ X(u), u \in [0,a];\ X(1) \}
  \mbox{\quad and \quad}
  {\cal H}:=
  \sigma \{ X(u), u \in [a,1] \}
  \ .
\]
  We have trivially 
  $\ds{\cal G}  
  \ci_{{\cal H}}
  {\cal H}$, 
  but the conditioning ${\cal H}$ cannot be replaced by the smaller 
  $\sigma$-field ${\cal H^\prime}:=\sigma \{ X(u), u \in [a,1]-[s,t]\}$.
  Indeed, it is easy to see that 
\[
  T:= \E\Big[\int_0^a X(u)\,du\ |\ {\cal H}\Big]
  =
  -\int_a^1 X(u)\,du\ ,
\]
  and one can check that $T$ is not ${\cal H^\prime}$-measurable:
  Choose $\omega^1,\omega^2\in C_0([0,1];\rr)$ such that $\omega^1\equiv\omega^2$
  on $\mathopen]s,t\mathclose[^c$ and 
  $\omega^1<\omega^2$ on $\mathopen]s,t\mathclose[$.                     
  An easy computation gives 
\[
  T(\omega^2)-T(\omega^1)
  =
  a\int_s^t (\omega^2-\omega^1)(u)\,du
  >0\ .
\]
  Now, due to the continuity of $T$ as a functional on $C_0([0,1];\rr)$, 
  it is possible to find two open balls, centred at $\omega^1$ 
  and $\omega^2$,   
  such that their images through $T$ take values in two disjoint intervals of
  $\rr$. We conclude that ${\cal G}$ and ${\cal H}$ are not conditionally 
  independent given ${\cal H^\prime}$.
  \qed
  }
\end{re}


\begin{thebibliography}{99}

\bibitem{af1}
Alabert, A., Ferrante, M.:
A conditional independence property for the solution of a
linear stochastic differential equation with lateral conditions.   
In: Stochastic analysis and related topics VI
(Prog.Probab. 42)
Boston: Birkh\"auser, 159-173 (1998).

\bibitem{afn}
Alabert, A., Ferrante, M., Nualart, D.:
Markov field property of stochastic differential equations.
Ann. Probab. {\bf 23}, 1262-1288 (1995).

\bibitem{am1}
Alabert, A., Marmolejo, M.A.:
Reciprocal property for a class of anticipating stochastic
differential equations.
Markov Process. Related Fields {\bf 5}, 331-356 (1999).

\bibitem{an1}
Alabert, A., Nualart, D.:
A second order Stratonovich differential equation with
boundary conditions.
Stochastic Processes Appl. {\bf 68}, 21-47 (1997).

\bibitem{ber}
Bernstein, S.: Sur les liaisons entre les grandeurs al\'eatoires.
In: Proc. Int. Cong. of Math., 288-309, Zurich (1932).

\bibitem{ChowLaso}
Chow, S-N., Lasota, A.:
On boundary value problems for ordinary differential equations.
Journal of Differential Equations {\bf 14}, 326-337 (1973).

\bibitem{co1}
Conti, R.:
Probl\`emes lin\'eaires pour les
\'equations diff\'erentielles ordinaires.
Math. Nachr. {\bf 23}, 161-178 (1961).

\bibitem{fe1}
Fe\v{c}kan, M.:
On the continuous dependence of solutions of nonlinear
equations,
J. Math. Anal. Appl. {\bf 194},
578-596 (1995).

\bibitem{fn1}
Ferrante, M., Nualart, D.:
An example of a non-Markovian stochastic two-point
boundary value problem,
Bernoulli {\bf 3}, 371-386 (1997).

\bibitem{FKL}
Frezza, R., Krener, A., Levy, C.:
{Gaussian reciprocal processes and selfadjoint stochastic
             differential equations of second order},
{Stochastics Stochastics Rep.}, {\bf 34}, 29-56 (1991).

\bibitem{gog}
Goggin, E.:
Convergence in distribution of conditional expectations.
Ann. Probab. {\bf 22}, 1097-1114 (1994).

\bibitem{hon}
H\"onig, C.S.:
The Green function of a linear differential
equation with a lateral condition.
Bull. Am. Math. Soc. {\bf 79}, 587-593 (1973).

\bibitem{jam}
Jamison, B.:
Reciprocal processes: The stationary Gaussian case.
The Annals of Mathematical Statistics {\bf 41}, 1624-1630 (1970).

\bibitem{kn1}
Knudsen, T.S.:
Convergence in the mean of conditional expectations.
(Preprint) 1998.

\bibitem{kre}
Krener, A.:
Reciprocal diffusions in flat space.
Probab. Theory Relat. Fields {\bf 107}, 243-281 (1997).


\bibitem{np1}
Nualart, D., Pardoux, E.:
Boundary value problems for stochastic differential equations.
Ann. Probab. {\bf 19}, 1118-1144 (1991).

\bibitem{np2}
Nualart, D., Pardoux, E.:
Second order stochastic differential equations with
Dirichlet boundary conditions.
Stochastic Processes Appl. {\bf 39}, 1-24 (1991).

\bibitem{op1}
Ocone, D., Pardoux, E.:
Linear stochastic differential equations with
boundary conditions.
Probab. Theory Relat. Fields {\bf 82}, 489-526 (1989).

\bibitem{rus}
Russek, A.:
Gaussian n-Markovian processes and stochastic
boundary value problems.
Z. Wahrscheinlichkeitstheor. Verw. Gebiete {\bf 53}, 117-122 (1980).

\bibitem{thi}
Thieullen, M.:
{Second order stochastic differential equations and
 non-{G}aussian reciprocal diffusions}.
{Probab. Theory Relat. Fields},
{\bf 97}, 231-257 (1993).

\bibitem{zam}
Zambrini, J.-C.:
{Probability and analysis in quantum physics}.
In: {Stochastic analysis, path integration and dynamics (Warwick, 1987)},
223-242, Longman Sci. Tech. (1989).

\end{thebibliography}
\end{document}